\documentclass[a4paper,11pt,reqno]{amsart}

\usepackage{amsfonts}
\usepackage{amsmath}
\usepackage{amssymb}

\usepackage{mathrsfs}
\usepackage[colorlinks]{hyperref}

\numberwithin{equation}{section}

\usepackage{graphicx}

\newtheorem{theorem}{Theorem}[section]
\newtheorem{defi}[theorem]{Definition}

\newtheorem{prop}[theorem]{Proposition}
\newtheorem{lemma}[theorem]{Lemma}

\newenvironment{pff}{\hspace*{-\parindent}{\bf Proof \,}}
{\hfill $\Box$ \vspace*{0.2cm}}

\def\R2n{{\mathbb R}^{2n}}

\def\R2{{\mathbb R}^2}
\def\R2n{{\mathbb R}^{2n}}

\def\N0{{\mathbb N}_{0}}

\def\l2h{{\ell^2(\hbar\mathbb Z^n)}}

\begin{document}
	\title[ $M$-Ellipticity of Fredholm Pseudo-Differential Operators on $L^p(\mathbb{R}^n)$ and G\r{a}rding's Inequality]{ $M$-Ellipticity of Fredholm Pseudo-Differential Operators on $L^p(\mathbb{R}^n)$ and G\r{a}rding's Inequality }
	\author[Aparajita Dasgupta]{Aparajita Dasgupta}
	\address{
		Aparajita Dasgupta:
		\endgraf
		Department of Mathematics
		\endgraf
		Indian Institute of Technology, Delhi, Hauz Khas
		\endgraf
		New Delhi-110016 
		\endgraf
		India
		\endgraf
		{\it E-mail address} {\rm adasgupta@maths.iitd.ac.in}
	}
	\author[Lalit Mohan]{Lalit Mohan}
	\address{
		Lalit Mohan:
		\endgraf
		Department of Mathematics
		\endgraf
		Indian Institute of Technology, Delhi
		\endgraf
		India
		\endgraf
		{\it E-mail address} {\rm mohanlalit871@gmail.com}
	}
	\date{\today}
	\subjclass{Primary 47G30; Secondary 35S05}
	\keywords{ $M$-ellipticity, Fredholm Operators, G\r{a}rding's Inequality }
	\begin{abstract}
		In this paper, we study the $M$-ellipticity of Fredholm pseudo-differential operators associated with weighted symbols on $L^p(\mathbb{R}^n)$, $1 < p < \infty$. We also prove the G\r{a}rding's inequality for $M$-elliptic operators and the hybrid class of pseudo-differential operators, namely  SG$M$-elliptic operators.
	\end{abstract}
	\maketitle
	\tableofcontents
	\section{Introduction}
	A general symbolic calculus and corresponding definition of global ellipticity for partial differential operators and pseudo-differential operators on noncompact
	manifolds represent a relevant issue of the modern Mathematical Analysis. In
	the case of the Euclidean space $\mathbb{R}^n$,  one can refer to the pseudo-differential operators  corresponding to the H\"{o}rmander symbol class, $S^{m}_{\rho,\delta}$, $m\in \mathbb{R}$, $0\leq \rho, \delta\leq 1,$ \cite{Hor, MWbook}.  The other is given by the so-called SG classes.  In \cite{Yu}, they are also called pseudo-differential operators with exit behavior. SG pseudo-differential operators and related topics can be found in \cite{Bo, AD, AD&MW, Ni, NiRo} and the references therein. In general to prove $L^p$ boundedness of pseudo-differential  operators associated with the H\"{o}rmander symbols, $S^{0}_{1,0},$ the key ingredient is the Mikjlin-H\"{o}rmander theorem about Fourier multipliers. But unfortunately, $L^{p}$-boundedness theorem for $p\neq 2$ doesnot hold for operators in $S^{0}_{\rho,0},$ for $0<\rho<1.$ Taylor in \cite{Ta} introduced a
	suitable symbol subclass, $M^m_{\rho,0},$ of $S^{m}_{\rho,0}$ giving a continuous operator for every $1<p<\infty$ and $0<\rho\leq 1.$  Further to this, Garello and Morando introduced the subclass $M^{m}_{\rho,\Lambda}$ of $S^{m}_{\rho,\Lambda},$ which are just weighted version of the ones introduced by Taylor in \cite{GM, GM1,GA}
	and developed the symbolic calculus for the associated pseudo-differential operators with
	many applications to study the regularity of multi-quasi-elliptic operators. Many studies on various properties of $M$-elliptic pseudo-differential operators can be found in \cite{GM, MWMelliptic, AM, AL, GA} and \cite{MWspectral}. But till now not much has been done on the $M$-ellipticity of Fredholm pseudo-differential operators.  In this paper, we prove that  Fredholm pseudo differential operators with weighted symbol are $M$-elliptic pseudo-differential operators on $\mathbb{R}^n$. Finally we prove the G\r{a}rding's Inequality for $M-$elliptic operators. Let us start with a few historical notes and basic definitions.\\
	
	A positive function $\Lambda \in C^{\infty}(\mathbb{R}^n)$ is said to be a weight function with polynomial growth if there exists positive constants $\mu_0$, $\mu_1$, $C_0$ and $C_1$ be such that $\mu_0 \leq \mu_1$ and $C_0 \leq C_1$, for which we have following conditions:
	$$C_{0}(1+|\xi|)^{\mu_0} \leq \Lambda(\xi) \leq C_{1}(1+|\xi|)^{\mu_1},\hspace{.3cm}  \xi \in \mathbb{R}^n.$$
	Also, we assume that there exists a real constant $\mu$ such that $\mu \geq \mu_1$ and for all multi-indices $\alpha$,$\gamma$ $\in \mathbb{Z}^n$ with $\gamma_j \in \{0,1\}, j = 1,2,...,n$, we can find a positive constant say $C_{\alpha,\gamma}$ such that
	\begin{eqnarray}\label{baisc-def}
		|\xi^\gamma(\partial^{\alpha+\gamma}\Lambda)(\xi)| \leq C_{\alpha,\gamma}\Lambda(\xi)^{1-\frac{1}{\mu}|\alpha|},\hspace{.3cm}  \xi \in \mathbb{R}^n.
	\end{eqnarray}
	Let us take $m \in \mathbb{R}$ and let $\rho \in (0,\frac{1}{\mu}]$. Then we say a function $\sigma \in C^{\infty}(\mathbb{R}^n \times \mathbb{R}^n)$ is in $S_{\rho,\Lambda}^m$ class if for all multi-indices $\alpha$ and $\beta$ we can find a positive constant $C_{\alpha,\beta}$ such that
	$$|(\partial_{x}^\alpha \partial_{\xi}^\beta\sigma)(x,\xi)| \leq C_{\alpha,\beta} \Lambda(\xi)^{m-\rho|\beta|},\hspace{.3cm}  x,\xi \in \mathbb{R}^n.$$
	We call such $\sigma$ be a symbol of order $m$ and type $\rho$ with weight $\Lambda$. Also, we say a symbol $\sigma \in M_{\rho,\Lambda}^m$ if for all multi-indices $\gamma$ with $\gamma_j \in \{0,1\}, j = 1,2,...,n$, we have
	$$\xi^{\gamma}(\partial_{\xi}^{\gamma}\sigma)(x,\xi) \in S_{\rho,\Lambda}^m.$$
	Let $\sigma \in S_{\rho,\Lambda}^m$. Define pseudo-differential operator $T_{\sigma}$ associate with symbol $\sigma$ by
	$$(T_{\sigma}\varphi)(x) = (2\pi)^{-\frac{n}{2}} \int_{\mathbb{R}^n} e^{i x \cdot \xi} \sigma(x,\xi) \hat{\varphi}(\xi)\, d\xi,\hspace{.3cm}  x \in \mathbb{R}^n,$$
	where $\varphi$ is in Schwartz space $\mathcal{S}$ and
	$$\hat{\varphi}(\xi) = (2\pi)^{-\frac{n}{2}} \int_{\mathbb{R}^n} e^{-i x \cdot \xi} \varphi(x)\, dx,\hspace{.3cm}  \xi \in \mathbb{R}^n.$$\\
	Now we give a brief introduction to the properties of $M$-elliptic pseudo-differential operators. For this, we will start by defining a class of symbols $S_{\rho,\Lambda}^{m}$, where $m$,  $\in \mathbb{R}$ and  $\rho \in (0,\frac{1}{\mu}]$.\\
	We say a symbol $\sigma \in C^{\infty}(\mathbb{R}^n \times \mathbb{R}^n)$ is in $S_{\rho,\Lambda}^{m}$ if for all multi-indices $\alpha$  we can find a positive constant $C_{\alpha}$ such that
	$$|(\partial_{x}^\alpha \partial_{\xi}^\beta\sigma)(x,\xi)| \leq C_{\alpha} \Lambda(\xi)^{m-\rho|\alpha|} ,\hspace{.3cm}  x,\xi \in \mathbb{R}^n.$$
	Also, we say a symbol $\sigma \in M_{\rho,\Lambda}^{m}$ if for all multi-indices $\gamma$ with $\gamma_j \in \{0,1\}, j = 1,2,...,n$, we have
	$$\xi^{\gamma}(\partial_{\xi}^{\gamma}\sigma)(x,\xi) \in S_{\rho,\Lambda}^{m}.$$
	Let $\sigma \in S_{\rho,\Lambda}^{m}$. Define pseudo-differential operator $T_{\sigma}$ associated with symbol $\sigma$ by
	\begin{eqnarray}\label{pseudo-def}
		(T_{\sigma}\varphi)(x) = (2\pi)^{-\frac{n}{2}} \int_{\mathbb{R}^n} e^{i x \cdot \xi} \sigma(x,\xi) \hat{\varphi}(\xi)\, d\xi,\hspace{.3cm}  x \in \mathbb{R}^n,
	\end{eqnarray}
	where $\varphi$ is in Schwartz space $\mathcal{S}$ and
	$$\hat{\varphi}(\xi) = (2\pi)^{-\frac{n}{2}} \int_{\mathbb{R}^n} e^{-i x \cdot \xi} \varphi(x)\, dx,\hspace{.3cm}  \xi \in \mathbb{R}^n.$$
	Note that it can be easily shown that $T_\sigma:\mathcal{S} \rightarrow \mathcal{S}$ is a continuous linear mapping.\\
	The following results can be found in \cite{AD&MW}.
	\begin{theorem}\label{product}
		Let $\sigma \in M_{\rho,\Lambda}^{m}$ and $\tau \in M_{\rho,\Lambda}^{\mu}$. Then $T_\sigma T_\tau = T_\lambda$, where $\lambda \in M_{\rho,\Lambda}^{m+\mu}$ and
		$$\lambda \sim \sum_{\mu} \frac{(-i)^{|\mu|}}{\mu!} (\partial_{\xi}^{\mu}\sigma) (\partial_{x}^{\mu}\tau).$$
		Here the asymptotic expansion means that for every positive integer $M$, there exists a positive integer $N$ such that
		$$\lambda - \sum_{|\mu|<N} \frac{(-i)^{|\mu|}}{\mu!} (\partial_{\xi}^{\mu}\sigma) (\partial_{x}^{\mu}\tau) \in M_{\rho,\Lambda}^{m+\mu-\rho M}.$$
	\end{theorem}
	\begin{theorem}\label{adjoint}
		Let $\sigma \in M_{\rho,\Lambda}^{m}$. Then the formal adjoint $T_{\sigma}^{\ast}$ of  $T_{\sigma}$ is the  pseudo-differential operator  $T_{\tau}$, where $\tau \in M_{\rho,\Lambda}^{m}$ and 
		$$\tau \sim \sum_{\mu} \frac{(-i)^{|\mu|}}{\mu!} \partial_{\xi}^{\mu}\partial_{x}^{\mu}\overline{\sigma}.$$
		Here the asymptotic expansion means that for every positive integer $M$, there exists a positive integer $N$ such that
		$$\tau - \sum_{|\mu|<N} \frac{(-i)^{|\mu|}}{\mu!}\partial_{\xi}^{\mu}\partial_{x}^{\mu}\overline{\sigma} \in M_{\rho,\Lambda}^{m-\rho M}.$$
	\end{theorem}
	Now using formal adjoint, we can extend the definition of a pseudo-differential operator from the Schwartz
	space $\mathcal{S}$ to the space $\mathcal{S}^{\prime}$ of all tempered distributions. For this, let $\sigma \in M_{\rho,\Lambda}^{m}$. Then for all $u$ in $\mathcal{S}^{\prime}$, we define $T_{\sigma}u:\mathcal{S} \rightarrow \mathbb{C} $ such that
	\begin{eqnarray}\label{schwartz-dual}
		(T_{\sigma}u)(\varphi) = u(\overline{T_{\sigma}^{\ast}\overline{\varphi}}).
	\end{eqnarray}
	It is easy to check that $T_{\sigma}$ maps $\mathcal{S}^{\prime}$ into $\mathcal{S}^{\prime}$ continuously. In fact, we have the following theorem.
	\begin{theorem}\label{boundness}
		Let $\sigma \in M_{\rho,\Lambda}^{0}$. Then $T_\sigma:L^p(\mathbb{R}^n) \rightarrow L^p(\mathbb{R}^n)$ is a bounded linear operator for $1 < p < \infty$.
	\end{theorem}
	Proof of this theorem can be found in from Theorem 1.4 in \cite{MWMelliptic}.\\\\
	Let $\sigma \in M_{\rho,\Lambda}^{m}$, where $m \in \mathbb{R}$. Then $\sigma$ is said to be $M$-elliptic if there exists positive constants $C$ and $R$ such that
	$$|\sigma(x,\xi)| \geq C\Lambda(\xi)^{m},\hspace{.3cm}   |\xi| \geq R.$$
	\begin{theorem}\label{ellip}
		Let $\sigma \in M_{\rho,\Lambda}^{m}$ be  $M$-elliptic. Then there exists a symbol $\tau \in M_{\rho,\Lambda}^{-m}$ such that
		$$T_\tau T_\sigma = I + R$$
		and
		$$T_\sigma T_\tau = I + S,$$
		where $R$ and $S$ are pseudo-differential operators with symbol in $\bigcap_{k \in \mathbb{R}} M_{\rho,\Lambda}^{k}$.
	\end{theorem}
	The pseudo-differential operator $T_{\tau}$ in Theorem \ref{ellip} is known as parametrix of the $M$-elliptic pseudo-differential operator $T_{\sigma}$.  \\
To make this paper self contained we recall  here the definition of 	the Sobolev spaces.\\	
		
For $m \in \mathbb{R}$ and $1\,<\,p\,<\,\infty$, we define the Sobolev space $H_{\Lambda}^{m,p}$ by $$H_{\Lambda}^{m,p} = \{u \in \mathcal{S}^{\prime} : J_{-m2}\,u \in \,L^p(\mathbb{R}^n)\}.$$ 
	Then $H_{\Lambda}^{m,p}$ is a Banach space in which the norm $\|\,\|_{m,p}$ is given by
	$$\|u\|_{m,p,\Lambda} = \|J_{-m}\,u\|_{L^p(\mathbb{R}^n)}, \hspace{0.3cm} u\,\in\,H_{\Lambda}^{m,p}.$$
	Note that $H_{\Lambda}^{0,p} = L^p(\mathbb{R}^n)$.\\
	
	All the above results and definitions for SG-$M$ elliptic case can be found in \cite{Al2}.\\

	The main aim of this paper is to prove the ellipticity of fredholm pseudo-differential operators associated with weighted symbols and to prove the G\r{a}rding's Inequality for $M$-elliptic operators.  In Section 2, first, we prove that Fredholm pseudo-differential operators are $M$-elliptic and then show that ellipticity of fredholmn pseudo-differential operators for the  hybrid case. In Section 3, we prove the G\r{a}rding's inequality for both $M$-elliptic operators and SG $M$-elliptic operators.

		\section{$M$-Ellipticity of Fredholm Pseudo-Differential Operators}
	In this section, our aim is to show that Fredholm pseudo-differential operator is $M$-elliptic
	on $L^p(\mathbb{R}^n)$. For this, we need some technical preparations which are done in \cite{AD} and \cite{MWbook}. 
	\begin{defi}\label{definition for R}
		Let $\lambda>0, \tau \geq 0$ and $x_{0}, \xi_{0} \in \mathbb{R}^{n}$. For $1<p<\infty$, we define the operator $R_{\lambda, \tau}\left(x_{0}, \xi_{0}\right): L^{p}\left(\mathbb{R}^{n}\right) \rightarrow L^{p}\left(\mathbb{R}^{n}\right)$, by
		$$
		\left(R_{\lambda, \tau}\left(x_{0}, \xi_{0}\right) u\right)(x)=\lambda^{\tau n / p}\, e^{i \lambda x  \cdot \xi_{0}}\, u\left(\lambda^{\tau}\left(x-x_{0}\right)\right), \quad x \in \mathbb{R}^{n},
		$$
		for all $u \in L^{p}\left(\mathbb{R}^{n}\right)$.
	\end{defi}
	\begin{prop}\label{surjectivity of R}
		The operator $R_{\lambda, \tau}\left(x_{0}, \xi_{0}\right): L^{p}\left(\mathbb{R}^{n}\right) \rightarrow L^{p}\left(\mathbb{R}^{n}\right)$ is a surjective isometry and the inverse is given by
		$$
		\left(R_{\lambda, \tau}\left(x_{0}, \xi_{0}\right)^{-1} u\right)(x)=\lambda^{-\tau n / p}\, e^{-i \lambda(x_{0}+\lambda^{-\top} x) \cdot \xi_{0}}\, u\left(x_{0}+\lambda^{-\tau} x\right), \quad x \in \mathbb{R}^{n},
		$$
		for all $u \in L^{p}\left(\mathbb{R}^{n}\right)$.
	\end{prop}
	\begin{prop}\label{weak convergence for R}
		Let $1<p<\infty$ and $\tau > 0$. Then, for all $u \in L^{p}\left(\mathbb{R}^{n}\right)$ and $v \in L^{p^{\prime}}\left(\mathbb{R}^{n}\right)$, where $\frac{1}{p}+\frac{1}{p^{\prime}}=1$,
		$$
		\left(R_{\lambda, \tau}\left(x_{0}, \xi_{0}\right) u, v\right) \rightarrow 0
		$$
		as $\lambda \rightarrow \infty$.
	\end{prop}
	The following theorem is one of the main theorems of this paper.
	\begin{theorem}\label{elliptic for 0}
		Let $\sigma \in M_{\rho,\Lambda}^{0}$ be such that $T_{\sigma}: L^{p}\left(\mathbb{R}^{n}\right) \rightarrow L^{p}\left(\mathbb{R}^{n}\right)$ is a Fredholm operator for $1<p<\infty .$ Then $\sigma$ is $M$-elliptic.
	\end{theorem}
	To prove this theorem, we need following lemma.
	\begin{lemma}\label{pdo with new symbol}
		Let $\sigma \in M_{\rho,\Lambda}^{m}$, where $m \in (-\infty,\infty).$ Then
		\begin{eqnarray}\label{lemma-eqn}
			R_{\lambda, \tau}\left(x_{0}, \xi_{0}\right)^{-1} T_{\sigma} R_{\lambda, \tau}\left(x_{0}, \xi_{0}\right)=T_{\sigma_{\lambda} \tau},
		\end{eqnarray}
		where
		\begin{eqnarray}\label{symbol-eqn}
			\sigma_{\lambda, \tau}(x, \eta)=\sigma\left(x_{0}+\lambda^{-\tau} x, \lambda \xi_{0}+\lambda^{\top} \eta\right), \quad x, \eta \in \mathbb{R}^{n}.
		\end{eqnarray}
		Moreover, if $\sigma \in M_{\rho,\Lambda}^{0}, \lambda \geq 1, 0 \leq \tau \leq \rho \mu_0/(1 + \rho \mu_0)$ and $\xi_{0} \neq 0$, then for all multi-indices $\alpha$ and $\beta$, there exists a positive constant $C_{\beta}$ such that
		\begin{eqnarray}\label{estimate-eqn}
			\left|\left(\partial_{x}^{\alpha} \partial_{\eta}^{\beta} \sigma_{\lambda, \tau}\right)(x, \eta)\right| \leq C_{\beta} p_{\alpha, \beta}(\sigma) \frac{\left(\Lambda(\eta)\right)^{\rho |\beta|}}{\left|\xi_{0}\right|^{\rho \mu_0 |\beta|}} \lambda^{-\tau|\alpha|} \lambda^{-(\rho \mu_0-(1+\rho \mu_0)\tau)|\beta|}, \quad x, \eta \in \mathbb{R}^{n}.
		\end{eqnarray}
		Here $p_{\alpha,\beta}$ denotes the corresponding norm in $M_{\rho,\Lambda}^{0}.$
	\end{lemma}
	\begin{pff}
		Fix $x_0,\xi_0 \in \mathbb{R}^{n}$. Let $u \in \mathcal{S}$. Then for all $x \in \mathbb{R}^{n}$, we can write $\left(R_{\lambda, \tau}\left(x_{0}, \xi_{0}\right)^{-1} T_{\sigma} R_{\lambda, \tau}\left(x_{0}, \xi_{0}\right) u\right)(x)$ in the form
		$$= e^{-i \lambda\left(x_{0}+\lambda^{-\top} x\right) \cdot \xi_{0}}\,(2 \pi)^{-n} \int_{\mathbb{R}^{n}} \int_{\mathbb{R}^{n}} e^{i \lambda\left(x_{0}+\lambda^{-\tau} x-y\right) \cdot \xi} \sigma\left(x_{0}+\lambda^{-\tau} x, \xi\right) e^{i \lambda y \cdot \xi_{0}} u\left(\lambda^{\tau}\left(y-x_{0}\right)\right) d y d \xi.$$
		Then by substituting $z = \lambda^{\tau} (y-x_0)$, the above expression takes the form
		$$= \lambda^{-\tau n}(2 \pi)^{-n} \int_{\mathbb{R}^{n}} \int_{\mathbb{R}^{n}} e^{i \lambda^{-\tau}\left(\xi-\lambda \xi_{0}\right) \cdot(x-z)} \sigma\left(x_{0}+\lambda^{-\tau} x, \xi\right) u(z) d z d \xi.$$
		Now by substituting $\eta = \lambda^{-\tau} (\xi-\lambda\xi_0)$, the above expression takes the form
		$$=(2 \pi)^{-n} \int_{\mathbb{R}^{n}} \int_{\mathbb{R}^{n}} e^{i(x-z) \cdot \eta} \sigma\left(x_{0}+\lambda^{-\tau} x, \lambda \xi_{0}+\lambda^{\tau} \eta\right) u(z) d z d \eta.$$
		Thus, we get \eqref{lemma-eqn} and \eqref{symbol-eqn}, as asserted. Now, using \eqref{symbol-eqn}, the chain rule and Peetre's inequality,
		$$
		\begin{aligned}
			\left|\left(\partial_{x}^{\alpha} \partial_{\eta}^{\beta} \sigma_{\lambda, \tau}\right)(x, \eta)\right| &=\left|\left(\partial_{x}^{\alpha} \partial_{\eta}^{\beta} \sigma\right)\left(x_{0}+\lambda^{-\tau} x, \lambda \xi_{0}+\lambda^{\tau} \eta\right) \lambda^{-\tau|\alpha|} \lambda^{\tau|\beta|}\right| \\
			& \leq p_{\alpha, \beta}(\sigma)\left(\Lambda(\lambda \xi_{0}+\lambda^{\tau} \eta)\right)^{-\rho|\beta|} \lambda^{-\tau|\alpha|} \lambda^{\tau|\beta|} \\
			& \leq C_{\beta} p_{\alpha, \beta}(\sigma)\left\langle\lambda \xi_{0}+\lambda^{\tau} \eta\right\rangle^{-\rho \mu_0 |\beta|} \lambda^{-\tau|\alpha|} \lambda^{\tau|\beta|} \\
			& \leq C_{\beta} p_{\alpha, \beta}(\sigma)\left(\lambda \xi_{0}\right\rangle^{-\rho \mu_0 |\beta|}\left\langle\lambda^{\tau} \eta\right\rangle^{\rho \mu_0 |\beta|} \lambda^{-\tau|\alpha|} \lambda^{\tau|\beta|}\\
			& \leq C_{\beta} p_{\alpha, \beta}(\sigma)\left|\xi_{0}\right|^{-\rho \mu_0 |\beta|}\langle\eta\rangle^{\rho \mu_0 |\beta|} \lambda^{-(\rho \mu_0 - (1+\rho \mu_0)\tau)|\beta|} \lambda^{-\tau|\alpha|}\\
			& \leq C_{\beta}^{\prime} p_{\alpha, \beta}(\sigma)\left|\xi_{0}\right|^{-\rho \mu_0 |\beta|} \left(\Lambda(\eta)\right)^{\rho |\beta|} \lambda^{-\left(\rho \mu_0 - (1+\rho \mu_0)\tau\right)|\beta|} \lambda^{-\tau|\alpha|}
		\end{aligned}
		$$
		which completes proof of Equation \eqref{estimate-eqn}.
	\end{pff}\\
	\textbf{Proof of Theorem \ref{elliptic for 0}} 
	Since $T_{\sigma}$ is a Fredholm operator, so by Theorem 20.5 in \cite{MWbook}, we can find a non-zero bounded linear operator $S$ on $L^{p}\left(\mathbb{R}^{n}\right)$ and a compact operator $K$ on $L^{p}\left(\mathbb{R}^{n}\right)$ such that
	$$
	S T_{\sigma}=I+K.
	$$
	Let $M$ be the set of all points $\xi$ in $\mathbb{R}^{n}$ such that there exists a point $x$ in $\mathbb{R}^{n}$ for which
	$$
	|\sigma(x, \xi)| \leq \frac{1}{2\|S\|}.
	$$
	Now, if $M$ is bounded, then there exists a positive number $R$ such that
	$$
	|\xi| < R, \quad \xi \in M.
	$$
	Thus, for each point $\xi \in \mathbb{R}^{n}$ with $|\xi| \geq R$, we get, for all $x \in \mathbb{R}^{n}$,
	$$
	|\sigma(x, \xi)| \geq \frac{1}{2\|S\|},
	$$
	and this implies that that $\sigma$ is $M$-elliptic. So, suppose that $M$ is not bounded. Then there exists a sequence $\left\{\left(x_{k}, \xi_{k}\right)\right\}$ in $\mathbb{R}^{n} \times \mathbb{R}^{n}$ such that
	$$
	\left|\xi_{k}\right| \rightarrow \infty
	$$
	as $k \rightarrow \infty$ and
	$$
	\left|\sigma\left(x_{k}, \xi_{k}\right)\right| \leq \frac{1}{2\|S\|}, \quad k=1,2, \ldots
	$$
	Thus, there exists a subsequence of $\left\{\left(x_{k}, \xi_{k}\right)\right\}$, again denoted by $\left\{\left(x_{k}, \xi_{k}\right)\right\}$, such that
	$$
	\sigma\left(x_{k}, \xi_{k}\right) \rightarrow \sigma_{\infty}
	$$
	for some complex number $\sigma_{\infty}$ as $k \rightarrow \infty .$ Therefore
	\begin{eqnarray}\label{one side inequality}
		\left|\sigma_{\infty}\right| \leq \frac{1}{2\|S\|} .
	\end{eqnarray}
	For $k=1,2, \ldots$, let $\lambda_{k}=\left|\xi_{k}\right| $. Then, by Lemma \ref{pdo with new symbol}, we have
	$$
	R_{\lambda_{k}, \tau}\left(x_{k}, \frac{\xi_{k}}{\left|\xi_{k}\right|}\right)^{-1} T_{\sigma} R_{\lambda_{h}, \tau}\left(x_{k}, \frac{\xi_{k}}{\left|\xi_{k}\right|}\right)=T_{\sigma_{\lambda_{k}, \tau}},
	$$
	where
	$$
	\sigma_{\lambda_{k}, \tau}(x, \eta)=\sigma\left(x_{k}+\lambda_{k}^{-\tau} x, \xi_{k}+\lambda_{k}^{\tau} \eta\right), \quad x, \eta \in \mathbb{R}^{n}.
	$$
	Let $\alpha$ and $\beta$ be arbitrary multi-indices. Then, by Equation \eqref{estimate-eqn}, there exists a positive constant $C_{\beta}$ such that
	\begin{eqnarray}\label{estimate for lambda_k}
		\left|\left(\partial_{x}^{\alpha} \partial_{\eta}^{\beta} \sigma_{\lambda_k, \tau}\right)(x, \eta)\right| \leq C_{\beta} p_{\alpha, \beta}(\sigma) \left(\Lambda(\eta)\right)^{\rho |\beta|}  \lambda_k^{-\tau|\alpha|} \lambda_k^{-(\rho \mu_0-(1+\rho \mu_0)\tau)|\beta|}, \quad x, \eta \in \mathbb{R}^{n}.
	\end{eqnarray}
	For $k=1,2, \ldots$, we define $\sigma_{k}^{\infty}$ by
	$$
	\sigma_{k}^{\infty}=\sigma_{\lambda_{k}, \tau}(0,0)=\sigma\left(x_{k}, \xi_{k}\right).
	$$
	Then, by using Theorem 7.3 in \cite{MWbook} and the estimate \eqref{estimate for lambda_k}, we get
	\begin{equation}\label{estimate by taylor formula}
		\begin{aligned}[b]
			&\hspace*{0.48cm}\left|\sigma_{\lambda_{k}, \tau}(x, \eta)-\sigma_{k}^{\infty}\right| \\
			&=\left|\sigma_{\lambda_{k}, \tau}(x, \eta)-\sigma_{\lambda_{k}, \tau}(0,0)\right|\\
			&=\left|\sum_{|\gamma+\mu|=1} x^{\gamma} \eta^{\mu} \int_{0}^{1}\left(\partial_{x}^{\gamma} \partial_{\eta}^{\mu} \sigma_{\lambda_{k}, \tau}\right)(\theta x, \theta \eta) d \theta\right|\\
			&\leq \sum_{|\gamma+\mu|=1}|x|^{|\gamma|}|\eta|^{|\mu|} \int_{0}^{1} C_{\mu} p_{\gamma, \mu}(\sigma)\left(\Lambda(\eta)\right)^{\rho |\mu|}  \lambda_k^{-\tau|\gamma|} \lambda_k^{-(\rho \mu_0-(1+\rho \mu_0)\tau)|\mu|}d \theta \rightarrow 0 
		\end{aligned}
	\end{equation}
	uniformly for $(x, \eta)$ on every compact subset $K$ of $\mathbb{R}^{n} \times \mathbb{R}^{n}$ as $k \rightarrow \infty$. Let $u \in \mathcal{S}$. Then
	$$
	\left(T_{\sigma_{\lambda_{k}}, \tau} u\right)(x)-\sigma_{k}^{\infty} u(x)=(2 \pi)^{-n / 2} \int_{\mathbb{R}^{n}} e^{i x \cdot \eta}\left(\sigma_{\lambda_{k}, \tau}(x, \eta)-\sigma_{k}^{\infty}\right) \hat{u}(\eta) d \eta
	$$
	for all $x \in \mathbb{R}^{n}$. By \eqref{estimate by taylor formula}, the assumption that $\sigma \in M_{\rho,\Lambda}^{0}$ and Lebesgue's dominated convergence theorem,
	$$
	\left(T_{\sigma_{\lambda_{k}, \tau}} u\right)(x) \rightarrow \sigma_{k}^{\infty} u(x)
	$$
	for all $x \in \mathbb{R}^{n}$ as $k \rightarrow \infty$. Moreover, for all $l \in \mathbb{N}$, using \eqref{estimate for lambda_k} and an integration by parts, we can find a positive constant $C_{\mu}$ for each $\mu$ with $|\mu| \leq 2 l$ such that
	$$
	\begin{aligned}
		&\left|\langle x)^{2 l}\left(T_{\sigma_{\lambda_{k}, \tau}} u\right)(x)\right| \\
		=&\left|\langle x\rangle^{2 l}(2 \pi)^{-n / 2} \int_{\mathbb{R}^{n}} e^{i x \cdot \eta} \sigma_{\lambda_{k}, \tau}(x, \eta) \hat{u}(\eta) d \eta\right| \\
		=&\left|(2 \pi)^{-n / 2} \int_{\mathbb{R}^{n}}\left(\left(I-\Delta_{\xi}\right)^{l} e^{i x \cdot \xi}\right) \sigma_{\lambda_{k}, \tau}(x, \xi) \hat{u}(\xi) d \xi\right| \\
		\leq& (2 \pi)^{-n / 2} \int_{\mathbb{R}^{n}}\left|\sum_{|\mu| \leq 2 l} \frac{1}{\mu !}\left(P^{\mu}(D) \hat{u}\right)(\xi)\right| C_{\mu}(\Lambda(\xi))^{|\mu|} d \xi\\
		\leq&(2 \pi)^{-n / 2} \int_{\mathbb{R}^{n}}\sum_{|\mu| \leq 2 l} \frac{C_{\mu}^{\prime}}{\mu !} |\hat{u}(\xi)| \langle \xi \rangle^{2l-|\mu|} \langle \xi \rangle^{\mu_1 |\mu|} d \xi
	\end{aligned}
	$$
	for all $x \in \mathbb{R}^{n}$, where $P(D)=(I-\Delta)^{l}$. So, there exists a positive constant $C$ such that
	$$
	\left|\left(T_{\sigma_{\lambda_{k}, \tau}} u\right)(x)\right|^{p} \leq C\langle x\rangle^{-2 l p}, \quad x \in \mathbb{R}^{n}.
	$$
	Now, if $2 l p>n$, then $\langle x\rangle^{-2 l p} \in L^{1}\left(\mathbb{R}^{n}\right)$. So, there exists a positive constant $C_{1}$ such that
	$$
	\left|\left(T_{\sigma_{\lambda_{k}, \tau}}-\sigma_{k}^{\infty}\right) u(x)\right|^p  \leq C_{1}\langle x\rangle^{-2 l p}, \quad x \in \mathbb{R}^{n}.
	$$
	Thus,
	$$
	T_{\sigma_{\lambda_{k}, \tau}} u \rightarrow \sigma_{\infty} u
	$$
	in $L^{p}\left(\mathbb{R}^{n}\right)$ as $k \rightarrow \infty$. Let $u$ be a nonzero function in $L^{p}\left(\mathbb{R}^{n}\right) .$ Since $R_{\lambda_{k}, \tau}\left(x_{k}, \frac{\xi_{k}}{\xi_{k} \mid}\right)$ is an isometry, it follows that
	\begin{equation}\label{estimate for u}
		\begin{aligned}[b]
			0<\|u\|_{p}=&\left\|R_{\lambda_{k}, \tau}\left(x_{k}, \frac{\xi_{k}}{\left|\xi_{k}\right|}\right) u\right\|_{p} \\
			=&\left\|\left(S T_{\sigma}-K\right) R_{\lambda_{k}, \tau}\left(x_{k}, \frac{\xi_{k}}{\left|\xi_{k}\right|}\right) u\right\|_{p} \\
			\leq &\left\|S T_{\sigma} R_{\lambda_{k}, \tau}\left(x_{k}, \frac{\xi_{k}}{\left|\xi_{k}\right|}\right) u\right\|{p}+\left\|{K R_{\lambda_{k}, \tau}}\left(x_{k}, \frac{\xi_{k}}{\left|\xi_{k}\right|}\right) u\right\|_{p} \\
			\leq&\left\|R_{\lambda_{k}, \tau}\left(x_{k}, \frac{\xi_{k}}{\left|\xi_{k}\right|}\right)^{-1} T_{\sigma} R_{\lambda_{k}, \tau}\left(x_{k}, \frac{\xi_{k}}{\left|\xi_{k}\right|}\right) u\right\|_{p}\|S\|
			+\left\|K R_{\lambda_{k}, \tau}\left(x_{k}, \frac{\xi_{k}}{\left|\xi_{k}\right|}\right) u\right\|_{p}.
		\end{aligned}
	\end{equation}
	Now, using the fact that $K$ is a compact operator and Proposition \ref{weak convergence for R}, it follows that
	$$
	\left\|K R_{\lambda_{k}, \tau}\left(x_{k}, \frac{\xi_{k}}{\left|\xi_{k}\right|}\right) u\right\|_{p} \rightarrow 0
	$$
	as $k \rightarrow \infty$. Then, by Equation \eqref{estimate for u},
	\begin{eqnarray}\label{right side inequality}
		\|u\|_{p} \leq\|S\|\left|\sigma_{\infty}\right|\|u\|_{p}.
	\end{eqnarray}
	Thus, Equation \eqref{one side inequality} and Equation \eqref{right side inequality} give the contradiction that
	$$
	\frac{1}{\|S\|} \leq\left|\sigma_{\infty}\right| \leq \frac{1}{2\|S\|},
	$$
	which completes the proof of the theorem.\\\\
	The preceding theorem can be generalized to the following theorem. 
	\begin{theorem}\label{elliptic for m}
		Let $\sigma \in M_{\rho,\Lambda}^{m},$ where $m \in (-\infty,\infty)$, and let $T_{\sigma}: H_{\Lambda}^{s, p} \rightarrow H_{\Lambda}^{s-m, p}$ be a Fredholm operator for some $s \in(-\infty, \infty) .$ Then $T_{\sigma}$ is an $M$-elliptic operator.
	\end{theorem}
	\begin{pff}
		By Theorem 1.6 in \cite{MWMelliptic}, the operators $T_{\sigma}: H^{s, p} \rightarrow H_{\Lambda}^{s-m, p}, J_{-s}: H_{\Lambda}^{s, p} \rightarrow L^{p}\left(\mathbb{R}^{n}\right), J_{s}: L^{p}\left(\mathbb{R}^{n}\right) \rightarrow H_{\Lambda}^{s, p}, J_{m-s}: H^{s-m, p} \rightarrow L^{p}\left(\mathbb{R}^{n}\right)$ and $J_{s-m}: L^{p}\left(\mathbb{R}^{n}\right) \rightarrow H^{s-m, p} $ are bounded linear operators. Here $J_{m}$ is a pseudo-differential operator with symbol $\sigma_{m} \in M_{\rho,\Lambda}^{-m}$, where
		$$
		\sigma_{m}(\xi)=(\Lambda(\xi))^{-m}, \quad \xi \in \mathbb{R}^{n}.
		$$
		Let
		\begin{eqnarray}\label{composition with symbol}
			J_{m-s} T_{\sigma} J_{s}=T_{\tau}.
		\end{eqnarray}
		Then
		$$
		T_{\tau}: L^{p}\left(\mathbb{R}^{n}\right) \rightarrow L^{p}\left(\mathbb{R}^{n}\right),
		$$
		where $\tau \in M_{\rho,\Lambda}^{0} .$ Since $J_{s}$ is bijective, it follows that $J_{s}$ is Fredholm and $M$-elliptic for all $s \in(-\infty, \infty) .$ So, by Theorem \ref{elliptic for 0}, $T_{\tau}$ is $M$-elliptic. By Equation \eqref{composition with symbol}, the fact that $J_{s}, s \in \mathbb{R}$, is bijective and product of two $M$-elliptic operators is again an $M$-elliptic operator, we get $T_{\sigma}$ is $M$-elliptic.
	\end{pff}\\
	The following theorem is a simple consequence of Theorem \ref{elliptic for 0} and Theorem \ref{elliptic for m}, which proves the $M$-ellipticity
	of the Fredholm SG pseudo-differential operators on $L^p(\mathbb{R}^n)$. Details about the Sobolev spaces, $H_{\Lambda}^{s_{1}, s_{2}, p},$ can be found in \cite{Al2}.
	\begin{theorem}
		Let $\sigma \in M_{\rho,\Lambda}^{m_{1}, m_{2}}$, where $m_{1}, m_{2} \in(-\infty, \infty)$ and let
		$$
		T_{\sigma}: H_{\Lambda}^{s_{1}, s_{2}, p} \rightarrow H_{\Lambda }^{s_{1}-m_{1}, s_{2}-m_{2}, p}
		$$
		is a Fredholm operator for some $s_1,s_2 \in (-\infty, \infty)$. Then $T_{\sigma}$ is  a SG $M$-elliptic operator. 
	\end{theorem}
	\section{G\r{a}rding's Inequality for $M$-elliptic operators}
	We begin with a definition.
	\begin{defi}
		Let $\sigma \in M_{\rho,\Lambda}^m$, where $m \in \mathbb{R}$. Then $\sigma$ is said to be strongly $M$- elliptic if there exist positive constants $C$ and $R$ for which
		$$\operatorname{Re}(\sigma(x,\xi)) \geq C\,\Lambda(\xi)^{m},\hspace{0.3cm} |\xi| \geq R.$$
	\end{defi}
	\begin{theorem}\textbf{(G\r{a}rding's Inequality For M-elliptic Operators)}\label{Gardings theorem}\\
		Let $\sigma \in M_{\rho,\Lambda}^{2m}$, where $m \in \mathbb{R}$, be strongly $M$- elliptic symbol. Then we can find positive constants $C^{\prime}$ and $C_s$ for every real number $s \geq \frac{\rho}{2}$ such that
		$$\operatorname{Re}(T_{\sigma}\,\varphi,\varphi) \geq C^{\prime}\,\|\varphi\|_{m,2,\Lambda}^2 - C_s\,\|\varphi\|_{m-s,2,\Lambda}^2,\hspace{0.3cm} \varphi \in \mathcal{S}.$$
	\end{theorem}
	To prove this theorem, we need the following three lemmas.
	\begin{lemma}\label{Compostion in S-class}
		Let $\sigma \in S_{\rho,\Lambda}^0$ and $F$ be a $C^{\infty}$- function on the complex plane $\mathbb{C}$. Then $F \circ \sigma \in S_{\rho,\Lambda}^0$.
	\end{lemma}
	\begin{pff}
		We need to prove that for any two multi-indices $\alpha$ and $\beta$, we can find a positive constant $C_{\alpha,\beta}$ such that
		\begin{eqnarray}\label{eqn to proof}
			|(\partial_{x}^\alpha \partial_{\xi}^\beta ( F \circ \sigma))(x,\xi)| \leq C_{\alpha,\beta}\,\Lambda(\xi)^{-\rho|\beta|},\hspace{0.3cm} x,\xi \in \mathbb{R}^n.
		\end{eqnarray}
		Let $\alpha$ and $\beta$ be two multi-indices such that $|\alpha + \beta| = 0$, i.e., $|\alpha| = 0$ and $|\beta| = 0$. Since $\sigma \in S_{\rho,\Lambda}^0$, so we can find a positive constant $C$ such that
		$$|\sigma(x,\xi)| \leq C, \hspace{0.3cm} x,\xi \in \mathbb{R}^n.$$
		Thus, $F \circ \sigma$ is a bounded and $C^\infty$ function on $\mathbb{R}^n \times \mathbb{R}^n$. So we can find another positive constant $C_0$ such that
		$$|(F \circ \sigma)(x,\xi)| \leq C_0, \hspace{0.3cm} x,\xi \in \mathbb{R}^n.$$
		Hence Equation \eqref{eqn to proof} is true for all multi-indices $\alpha$ and $\beta$ with $|\alpha + \beta| = 0$.
		Now, suppose that Equation \eqref{eqn to proof} is true for all $C^\infty$ functions $F$ on $\mathbb{C}$, $\sigma \in S_{\rho,\Lambda}^0$ and multi-indices $\alpha$ and $\beta$ with $|\alpha + \beta| \leq k$. Let $\alpha$ and $\beta$ be multi-indices with $|\alpha + \beta| = k + 1$. We first suppose that
		$$\partial_{x}^{\alpha}\,\partial_{\xi}^{\beta} = \partial_{x}^{\alpha}\,\partial_{\xi}^{\gamma}\,\partial_{\xi_j}$$
		for some multi-index $\gamma$ with $|\alpha + \gamma| = k$ and for some $j = 1,2,...,n$. Then, by chain rule,
		$$(\partial_{x}^{\alpha}\,\partial_{\xi}^{\beta}(F \circ \sigma))(x,\xi) = \partial_{x}^{\alpha}\,\partial_{\xi}^{\gamma}\{(F_1 \circ \sigma)\,\partial_{\xi_j}\sigma + (F_2 \circ \sigma)\,\partial_{\xi_j}\sigma\}(x,\xi)$$
		for all $x$ and $\xi$ in $\mathbb{R}^n$, where $F_1$ and $F_2$ are the partial derivatives of $F$ with respect to the first and second variables respectively. Now, by Leibniz’s formula and the induction hypothesis, we can find two positive constants $C_{\delta_0,\delta}$ and $C_{\alpha,\delta_0,\gamma,\delta,j}$ such that
		\begin{align*}
			&|(\partial_{x}^{\alpha}\,\partial_{\xi}^{\gamma}\{(F_1 \circ \sigma)\,\partial_{\xi_j}\sigma\})(x,\xi)|\\
			&\leq \sum_{\substack{\delta_0 \leq \alpha\\ \delta \leq \gamma}} {\alpha \choose \delta_0}\, {\gamma \choose \delta}\, |(\partial_{x}^{\delta_0}\,\partial_{\xi}^{\delta}(F_1 \circ \sigma))(x,\xi)|\,|(\partial_{x}^{\alpha - \delta_0}\,\partial_{\xi}^{\gamma - \delta}\,(\partial_{\xi_j}\sigma))(x,\xi)|\\
			&\leq \sum_{\substack{\delta_0 \leq \alpha\\ \delta \leq \gamma}} {\alpha \choose \delta_0}\, {\gamma \choose \delta}\, C_{\delta_0,\delta}\,\Lambda(\xi)^{-\rho|\delta|}\,C_{\alpha,\delta_0,\gamma,\delta,j}\,\Lambda(\xi)^{-\rho(|\gamma| - |\delta| + 1)}\\
			&\leq C_{\alpha,\gamma,j}\,\Lambda(\xi)^{-\rho(|\gamma| + 1)}\\
			&= C_{\alpha,\gamma,j}\,\Lambda(\xi)^{-\rho(|\beta|)}, \hspace{0.3cm} x,\xi \in \mathbb{R}^n,
		\end{align*}
		where $$C_{\alpha,\gamma,j} = \sum_{\substack{\delta_0 \leq \alpha\\ \delta \leq \gamma}} {\alpha \choose \delta_0}\, {\gamma \choose \delta}\,C_{\delta_0,\delta}\,C_{\alpha,\delta_0,\gamma,\delta,j}.$$
		Similarly, we can find another positive constant $C_{\alpha,\gamma,j}^{\prime}$ such that
		$$|(\partial_{x}^{\alpha}\,\partial_{\xi}^{\gamma}\{(F_2 \circ \sigma)\,\partial_{\xi_j}\sigma\})(x,\xi)| \leq C_{\alpha,\gamma,j}^{\prime}\,\Lambda(\xi)^{-\rho(|\beta|)}, \hspace{0.3cm} x,\xi \in \mathbb{R}^n.$$
		Hence  $$|(\partial_{x}^{\alpha}\,\partial_{\xi}^{\beta}(F \circ \sigma))(x,\xi)| \leq (C_{\alpha,\gamma,j} + C_{\alpha,\gamma,j}^{\prime})\,\Lambda(\xi)^{-\rho(|\beta|)}, \hspace{0.3cm} x,\xi \in \mathbb{R}^n.$$
		Now, we assume that $$\partial_{x}^{\alpha}\,\partial_{\xi}^{\beta} = \partial_{x}^{\gamma}\,\partial_{x_j}\,\partial_{\xi}^{\beta}$$
		for some multi-index $\gamma$ with $|\gamma + \beta| = k$ and for some $j = 1,2,...,n$. Then, as before, we can find a positive constant $C_{\alpha,\gamma,j}^{\prime\prime}$ such that 
		$$|(\partial_{x}^{\alpha}\,\partial_{\xi}^{\beta}(F \circ \sigma))(x,\xi)| \leq C_{\alpha,\gamma,j}^{\prime\prime}\,\Lambda(\xi)^{-\rho(|\beta|)}, \hspace{0.3cm} x,\xi \in \mathbb{R}^n.$$
		So, by the principle of mathematical induction, Equation \eqref{eqn to proof} follows.
	\end{pff}
	\begin{lemma}\label{Compostion in M-class}
		Let $\sigma \in M_{\rho,\Lambda}^0$ and $F$ be a $C^{\infty}$- function on the complex plane $\mathbb{C}$. Then $F \circ \sigma \in M_{\rho,\Lambda}^0$.
	\end{lemma}
	\begin{pff}
		Let $S = \{\gamma \in \mathbb{Z}_{0}^{n} : \gamma_j \in \{0,1\}, j = 1,2,...,n; n \in \mathbb{N}\}$, where $\mathbb{Z}_0 = \mathbb{Z}_{+}\,\cup\,\{0\}.$\\
		Then, we need to show that, for all $\gamma \in S$, we have $$\xi^\gamma\,\partial_{\xi}^\gamma\,(F \circ \sigma)(x,\xi) \in 	S_{\rho,\Lambda}^{0},$$
		i.e., for all multi-indices $\alpha,\beta$, we can find a positive constant $C_{\alpha,\beta}$ such that
		\begin{eqnarray}\label{M-class proof}
			|\partial_{x}^\alpha\,\partial_{\xi}^\beta(\xi^\gamma\,\partial_{\xi}^\gamma\,(F \circ \sigma)(x,\xi))| \leq C_{\alpha,\beta}\,(\Lambda(\xi))^{-\rho\,|\beta|}, \hspace{0.3cm} \forall \gamma \in S, \hspace{0.3cm} x,\xi \in \mathbb{R}^n.
		\end{eqnarray}
		Let $|\alpha + \beta| = 0.$ This implies $|\alpha| = 0$ and $|\beta| = 0.$ Then we need to show that there exists a positive constant $C_{00}$ such that
		\begin{eqnarray}\label{M-class for 0}
			|(\xi^\gamma\,\partial_{\xi}^\gamma\,(F \circ \sigma)(x,\xi))| \leq C_{00}, \hspace{0.3cm} \forall \gamma \in S, \hspace{0.3cm} x,\xi \in \mathbb{R}^n.
		\end{eqnarray}
		Let $|\gamma| = 0$. Then by Lemma \ref{Compostion in S-class}, we can find a positive constant $C_{00}$ such that
		$$|(F \circ \sigma)(x,\xi)| \leq C_{00},\hspace{0.3cm} x,\xi \in \mathbb{R}^n.$$
		Now, suppose that Equation \eqref{M-class for 0} is true for all $C^\infty$ functions $F$ on $\mathbb{C}$, $\sigma \in M_{\rho,\Lambda}^0$ and multi-indices $\gamma$  with $|\gamma| \leq k$. Let $\gamma$ be a multi-index with $|\gamma| = k + 1$. We suppose that $\partial_{\xi}^\gamma = \partial_{\xi}^{\gamma_0}\,\partial_{\xi_j}$, for some $j = 1,2,...,n$ and multi-index $\gamma_0$ with $|\gamma_0| = k$. Then by chain rule,
		$${\xi}^{\gamma}\,\partial_{\xi}^{\gamma}(F \circ \sigma)(x,\xi) = {\xi}^{\gamma}\,\partial_{\xi}^{\gamma_0}\{(F_1 \circ \sigma)\,\partial_{\xi_j}\sigma + (F_2 \circ \sigma)\,\partial_{\xi_j}\sigma\}(x,\xi)$$
		for all $x$ and $\xi$ in $\mathbb{R}^n$, where $F_1$ and $F_2$ are the partial derivatives of $F$ with respect to the first and second variables respectively. Now, by Leibniz’s formula,
		\begin{align*}
			&{\xi}^{\gamma}\,\partial_{\xi}^{\gamma_0}((F_1 \circ \sigma)(x,\xi)\,.\,\partial_{\xi_j}(\sigma(x,\xi)))\\
			&= {\xi}^{\gamma}\,\sum_{\delta_0 \leq \gamma_0} {\delta_0 \choose \gamma_0}\,\partial_{\xi}^{\delta_0}(F_1 \circ \sigma)(x,\xi)\,.\,\partial_{\xi}^{\gamma_0-\delta_0}(\partial_{\xi_j}\sigma(x,\xi))\\
			&= \sum_{\delta_0 \leq \gamma_0} {\delta_0 \choose \gamma_0}\,({\xi}^{\delta_0}\,\partial_{\xi}^{\delta_0}(F_1 \circ \sigma)(x,\xi))\,.\,({\xi}^{\gamma_0-\delta_0+1}\partial_{\xi}^{\gamma_0-\delta_0+1}(\sigma(x,\xi))).
		\end{align*}
		Thus by induction hypothesis and using the fact that $\sigma \in M_{\rho,\Lambda}^{0}$, we can find positive constants $C_{\delta_0}$ and $C_{\gamma_0,\delta_0}$ such that 
		\begin{align*}
			&|{\xi}^{\gamma}\,\partial_{\xi}^{\gamma_0}((F_1 \circ \sigma)(x,\xi)\,.\,\partial_{\xi_j}(\sigma(x,\xi)))|\\
			&\leq \sum_{\delta_0 \leq \gamma_0} {\delta_0 \choose \gamma_0}\,|({\xi}^{\delta_0}\,\partial_{\xi}^{\delta_0}(F_1 \circ \sigma)(x,\xi))|\,.\,|({\xi}^{\gamma_0-\delta_0+1}\partial_{\xi}^{\gamma_0-\delta_0+1}(\sigma(x,\xi)))|\\
			&\leq \sum_{\delta_0 \leq \gamma_0} {\delta_0 \choose \gamma_0}\,C_{\delta_0}\,.\,C_{\gamma_0,\delta_0}.
		\end{align*}
		So, we have
		$$|{\xi}^{\gamma}\,\partial_{\xi}^{\gamma_0}((F_1 \circ \sigma)(x,\xi)\,.\,\partial_{\xi_j}(\sigma(x,\xi)))| \leq C_{00}^{\prime},\hspace{0.3cm} x,\xi \in \mathbb{R}^n,$$
		where $$C_{00}^{\prime} = \sum_{\delta_0 \leq \gamma_0} {\delta_0 \choose \gamma_0}\,C_{\delta_0}\,.\,C_{\gamma_0,\delta_0}.$$ Similarly, we can find a positive constant $C_{00}^{\prime\prime}$ such that
		$$|{\xi}^{\gamma}\,\partial_{\xi}^{\gamma_0}((F_2 \circ \sigma)(x,\xi)\,.\,\partial_{\xi_j}(\sigma(x,\xi)))| \leq C_{00}^{\prime\prime},\hspace{0.3cm} x,\xi \in \mathbb{R}^n.$$
		Hence
		$$|{\xi}^{\gamma}\,\partial_{\xi}^{\gamma}(F \circ \sigma)(x,\xi)| \leq (C_{00}^{\prime} + C_{00}^{\prime\prime}),\hspace{0.3cm} x,\xi \in \mathbb{R}^n.$$
		Thus, by the principle of mathematical induction,
		Equation \eqref{M-class for 0} follows.\\
		Now, suppose that Equation \eqref{M-class proof} is true for all $C^\infty$ functions $F$ on $\mathbb{C}$, $\sigma \in M_{\rho,\Lambda}^0$ and multi-indices $\alpha$ and $\beta$ with $|\alpha + \beta| \leq k$. Let $\alpha$ and $\beta$ be multi-indices with $|\alpha + \beta| = k + 1$. Then we need to show that we can find a positive constant $C_{\alpha,\beta}$ such that
		\begin{eqnarray}\label{M-class for k+1}
			|\partial_{x}^\alpha\,\partial_{\xi}^\beta(\xi^\gamma\,\partial_{\xi}^\gamma\,(F \circ \sigma)(x,\xi))| \leq C_{\alpha,\beta}\,(\Lambda(\xi))^{-\rho\,|\beta|}, \hspace{0.3cm} \forall \gamma \in S, \hspace{0.3cm} x,\xi \in \mathbb{R}^n.
		\end{eqnarray}
		We first suppose that
		$$\partial_{x}^{\alpha}\,\partial_{\xi}^{\beta} = \partial_{x}^{\alpha}\,\partial_{\xi}^{\beta_0}\,\partial_{\xi_j}$$
		for some multi-index $\beta_0$ with $|\alpha + \beta_0| = k$ and for some $j = 1,2,...,n$.\\
		Let $|\gamma| = 0$. Then by Lemma \ref{Compostion in S-class}, we can find a positive constant $C_{\alpha,\beta}$ such that
		$$|(\partial_{x}^\alpha \partial_{\xi}^\beta ( F \circ \sigma))(x,\xi)| \leq C_{\alpha,\beta}\,\Lambda(\xi)^{-\rho|\beta|},\hspace{0.3cm} \xi \in \mathbb{R}^n.$$
		Now, suppose that Equation \eqref{M-class for k+1} is true for all $C^\infty$ functions $F$ on $\mathbb{C}$, $\sigma \in M_{\rho,\Lambda}^0$ and multi-indices $\gamma$  with $|\gamma| \leq k$. Let $\gamma$ be a multi-index with $|\gamma| = k + 1$. We suppose that $\partial_{\xi}^\gamma = \partial_{\xi}^{\gamma_0}\,\partial_{\xi_i}$, for some $i = 1,2,...,n$ and multi-index $\gamma_0$ with $|\gamma_0| = k$. Then, by chain rule,
		\begin{align*}
			\partial_{x}^\alpha\,\partial_{\xi}^\beta(\xi^\gamma\,\partial_{\xi}^\gamma\,(F \circ \sigma)(x,\xi))
			&=\partial_{x}^\alpha\,\partial_{\xi}^\beta(\xi^\gamma\,\partial_{\xi}^{\gamma_0}\,\partial_{\xi_i}(F \circ \sigma)(x,\xi)))\\
			&=\partial_{x}^{\alpha}\,\partial_{\xi}^{\beta}\,({\xi}^{\gamma}\,\partial_{\xi}^{\gamma_0}\{(F_1 \circ \sigma)\,\partial_{\xi_i}\sigma + (F_2 \circ \sigma)\,\partial_{\xi_i}\sigma\}(x,\xi)).
		\end{align*}
		Thus, by Leibniz’s formula,
		\begin{align*}
			&\partial_{x}^{\alpha}\,\partial_{\xi}^{\beta}\,({\xi}^{\gamma}\,\partial_{\xi}^{\gamma_0}((F_1 \circ \sigma)\,.\,\partial_{\xi_i}\sigma)(x,\xi))\\
			&= \partial_{x}^{\alpha}\,\partial_{\xi}^{\beta}\,\left({\xi}^{\gamma}\,\left(\sum_{\delta_0 \leq \gamma_0} {\delta_0 \choose \gamma_0}\,\partial_{\xi}^{\delta_0}(F_1 \circ \sigma)(x,\xi)\,.\,\partial_{\xi}^{\gamma_0-\delta_0}(\partial_{\xi_i}\sigma(x,\xi))\right)\right)\\
			&= \partial_{x}^{\alpha}\,\partial_{\xi}^{\beta_0}\,\partial_{\xi_j}\,\left(\sum_{\delta_0 \leq \gamma_0} {\delta_0 \choose \gamma_0}\,(\xi^{\delta_0}\,\partial_{\xi}^{\delta_0}(F_1 \circ \sigma)(x,\xi))\,.\,({\xi}^{\gamma_0-\delta_0+1}\partial_{\xi}^{\gamma_0-\delta_0+1}\sigma(x,\xi))\right)\\
			&= \partial_{x}^{\alpha}\,\partial_{\xi}^{\beta_0}\,\left(\sum_{\delta_0 \leq \gamma_0} {\delta_0 \choose \gamma_0}\,\partial_{\xi_j}\left((\xi^{\delta_0}\,\partial_{\xi}^{\delta_0}(F_1 \circ \sigma)(x,\xi))\,.\,({\xi}^{\gamma_0-\delta_0+1}\partial_{\xi}^{\gamma_0-\delta_0+1}\sigma(x,\xi))\right)\right)\\
			&= \partial_{x}^{\alpha}\,\partial_{\xi}^{\beta_0}\,\left(\sum_{\delta_0 \leq \gamma_0} {\delta_0 \choose \gamma_0}\,\partial_{\xi_j}\left(\xi^{\delta_0}\,\partial_{\xi}^{\delta_0}(F_1 \circ \sigma)(x,\xi)\right)\,.\,({\xi}^{\gamma_0-\delta_0+1}\partial_{\xi}^{\gamma_0-\delta_0+1}\sigma(x,\xi))\right)\\
			&\hspace{0.4cm} +  \partial_{x}^{\alpha}\,\partial_{\xi}^{\beta_0}\,\left(\sum_{\delta_0 \leq \gamma_0} {\delta_0 \choose \gamma_0}\,(\xi^{\delta_0}\,\partial_{\xi}^{\delta_0}(F_1 \circ \sigma)(x,\xi))\,.\,\partial_{\xi_j}\left(({\xi}^{\gamma_0-\delta_0+1}\partial_{\xi}^{\gamma_0-\delta_0+1}\sigma(x,\xi))\right)\right).
		\end{align*}
		First, consider
		\begin{align*}
			&\partial_{x}^{\alpha}\,\partial_{\xi}^{\beta_0}\,\left(\sum_{\delta_0 \leq \gamma_0} {\delta_0 \choose \gamma_0}\,\partial_{\xi_j}\left(\xi^{\delta_0}\,\partial_{\xi}^{\delta_0}(F_1 \circ \sigma)(x,\xi)\right)\,.\,({\xi}^{\gamma_0-\delta_0+1}\partial_{\xi}^{\gamma_0-\delta_0+1}\sigma(x,\xi))\right)\\
			&= \sum_{\delta_0 \leq \gamma_0}\sum_{p_0 \leq \beta_0}\sum_{p \leq \alpha} {\delta_0 \choose \gamma_0}{p_0 \choose \beta_0}{p \choose \alpha}\,\partial_{x}^{p}\,\partial_{\xi}^{p_0+1}\left(\xi^{\delta_0}\,\partial_{\xi}^{\delta_0}(F_1 \circ \sigma)(x,\xi)\right).\,\partial_{x}^{\alpha-p}\,\partial_{\xi}^{\beta_0-p_0}\left({\xi}^{\gamma_0-\delta_0+1}\partial_{\xi}^{\gamma_0-\delta_0+1}\sigma(x,\xi)\right).
		\end{align*}
		Thus by induction hypothesis and using the fact that $\sigma \in M_{\rho,\Lambda}^{0}$, we can find positive constants $C_{p,p_0}$ and $C_{p,p_0,\alpha,\beta_0}$ such that
		\begin{align*}
			&\left|\partial_{x}^{\alpha}\,\partial_{\xi}^{\beta_0}\,\left(\sum_{\delta_0 \leq \gamma_0} {\delta_0 \choose \gamma_0}\,\partial_{\xi_j}\left(\xi^{\delta_0}\,\partial_{\xi}^{\delta_0}(F_1 \circ \sigma)(x,\xi)\right)\,.\,({\xi}^{\gamma_0-\delta_0+1}\partial_{\xi}^{\gamma_0-\delta_0+1}\sigma(x,\xi))\right)\right|\\
			&\leq \sum_{\delta_0 \leq \gamma_0}\sum_{p_0 \leq \beta_0}\sum_{p \leq \alpha} {\delta_0 \choose \gamma_0}{p_0 \choose \beta_0}{p \choose \alpha}\,C_{p,p_0}\,(\Lambda(\xi))^{-\rho|p_0+1|}\,.\,C_{p,p_0,\alpha,\beta_0}\,(\Lambda(\xi))^{-\rho|\beta_0-p_0|}\\
			&= C_{\alpha,\beta}^{\prime}\,(\Lambda(\xi))^{-\rho|\beta|},\hspace{0.3cm} x,\xi \in \mathbb{R}^n,
		\end{align*}
		where $$C_{\alpha,\beta}^{\prime} = \sum_{\delta_0 \leq \gamma_0}\sum_{p_0 \leq \beta_0}\sum_{p \leq \alpha} {\delta_0 \choose \gamma_0}{p_0 \choose \beta_0}{p \choose \alpha}\,C_{p,p_0}\,C_{p,p_0,\alpha,\beta_0}.$$
		Similarly, we can find a positive constant $C_{\alpha,\beta}^{\prime\prime}$ such that
		\begin{align*}
			&\left|\partial_{x}^{\alpha}\,\partial_{\xi}^{\beta_0}\,\left(\sum_{\delta_0 \leq \gamma_0} {\delta_0 \choose \gamma_0}\,(\xi^{\delta_0}\,\partial_{\xi}^{\delta_0}(F_1 \circ \sigma)(x,\xi))\,.\,\partial_{\xi_j}\left(({\xi}^{\gamma_0-\delta_0+1}\partial_{\xi}^{\gamma_0-\delta_0+1}\sigma(x,\xi))\right)\right)\right|\\
			&\leq C_{\alpha,\beta}^{\prime\prime}\,(\Lambda(\xi))^{-\rho|\beta|},\hspace{0.3cm} x,\xi \in \mathbb{R}^n.
		\end{align*}
		Hence
		$$|\partial_{x}^{\alpha}\,\partial_{\xi}^{\beta}\,({\xi}^{\gamma}\,\partial_{\xi}^{\gamma_0}((F_1 \circ \sigma)\,.\,\partial_{\xi_i}\sigma)(x,\xi))| \leq (C_{\alpha,\beta}^{\prime} + C_{\alpha,\beta}^{\prime\prime})\,(\Lambda(\xi))^{-\rho|\beta|},\hspace{0.3cm} x,\xi \in \mathbb{R}^n.$$
		Similarly, we can find a positive constant $C_{\alpha,\beta}^{\prime\prime\prime}$ such that
		$$|\partial_{x}^{\alpha}\,\partial_{\xi}^{\beta}\,({\xi}^{\gamma}\,\partial_{\xi}^{\gamma_0}((F_2 \circ \sigma)\,.\,\partial_{\xi_i}\sigma)(x,\xi))| \leq C_{\alpha,\beta}^{\prime\prime\prime}\,(\Lambda(\xi))^{-\rho|\beta|},\hspace{0.3cm} x,\xi \in \mathbb{R}^n.$$
		Hence $$|\partial_{x}^\alpha\,\partial_{\xi}^\beta(\xi^\gamma\,\partial_{\xi}^\gamma\,(F \circ \sigma)(x,\xi))| \leq (C_{\alpha,\beta}^{\prime} + C_{\alpha,\beta}^{\prime\prime} + C_{\alpha,\beta}^{\prime\prime\prime})\,(\Lambda(\xi))^{-\rho|\beta|},\hspace{0.3cm} x,\xi \in \mathbb{R}^n.$$
		Now, we assume that $$\partial_{x}^{\alpha}\,\partial_{\xi}^{\beta} = \partial_{x}^{\gamma}\,\partial_{x_j}\,\partial_{\xi}^{\beta}$$
		for some multi-index $\gamma$ with $|\gamma + \beta| = k$ and for some $j = 1,2,...,n$. Then, as before, we can find a positive constant $C_{\alpha,\beta}^{\prime\prime\prime\prime}$ such that 
		$$|(\partial_{x}^{\alpha}\,\partial_{\xi}^{\beta}(F \circ \sigma))(x,\xi)| \leq C_{\alpha,\beta}^{\prime\prime\prime\prime}\,\Lambda(\xi)^{-\rho(|\beta|)}, \hspace{0.3cm} x,\xi \in \mathbb{R}^n.$$
		So, by the principle of mathematical induction, Equation \eqref{M-class for k+1} follows.\\
		This completes the proof of the Equation \eqref{M-class proof}.
	\end{pff}
	\begin{lemma}\label{strong elliptic inequality}
		Let $\sigma$ be a strongly $M$-elliptic symbol in $M_{\rho,\Lambda}^{2m}$, where $m \in \mathbb{R}$. Then we can find two positive constants $\eta$ and $\kappa$ such that $$\operatorname{Re}(\sigma (x,\xi)) \geq \eta\,(\Lambda(\xi))^{2m} - \kappa\,(\Lambda(\xi))^{2m-\rho},\hspace{0.3cm} x,\xi \in \mathbb{R}^n.$$
	\end{lemma}
	\begin{pff}
		By strong ellipticity, there exist positive constants $C$ and $R$ such that $$\operatorname{Re}(\sigma (x,\xi)) \geq C\,(\Lambda(\xi))^{2m}, \hspace{0.3cm} |\xi | \geq R.$$
		Since $\sigma \in M_{\rho,\Lambda}^{2m}$, we can find a positive constant $K$ such that $$|\sigma(x,\xi)| \leq K\,(\Lambda(\xi))^{2m}, \hspace{0.3cm} x, \xi \in R^n.$$
		Therefore, if $m \geq 0$, then \begin{eqnarray}\label{for positive m}
			|\operatorname{Re}(\sigma (x,\xi))| \leq K\,(\Lambda(\xi))^{2m} \leq K_{2m}\,(1 + |\xi|)^{2 \mu_1. m} \leq K_{2m}\,(1 + R)^{2 \mu_1. m},\hspace{0.3cm} |\xi| \leq R,
		\end{eqnarray}
		where $$K_{2m} = K\,.\,C_{1}^{2m}.$$
		And if $m < 0$, then
		\begin{eqnarray}\label{for negative m}
			|\operatorname{Re}(\sigma (x,\xi))| \leq K\,(\Lambda(\xi))^{2m} \leq K,\hspace{0.3cm} |\xi| \leq R.
		\end{eqnarray}
		By Equation \eqref{for positive m} and Equation \eqref{for negative m}, for given $m \in \mathbb{R}$, we can find a positive constant $M$ such that $$\operatorname{Re}(\sigma (x,\xi)) \geq -M, \hspace{0.3cm} |\xi| \leq R.$$
		Since $\frac{\operatorname{Re}\sigma}{(\Lambda(.))^{2m-\rho}}$ is continuous on the compact set $\{\xi \in \mathbb{R}^n : |\xi| \leq R \}$, so we can find a positive constant $\kappa$ such that $$\frac{\operatorname{Re}\sigma}{(\Lambda(\xi))^{2m-\rho}} > - \kappa, \hspace{0.3cm} |\xi| \leq R.$$
		Therefore $$\operatorname{Re}(\sigma (x,\xi)) + \kappa\,(\Lambda(\xi))^{2m-\rho} > 0, \hspace{0.3cm}|\xi| \leq R.$$ 
		Since $\frac{\operatorname{Re}\sigma + \kappa\,(\Lambda(.))^{2m-\rho}}{(\Lambda(.))^{2m}}$ is a positive and continuous function on the compact set $\{\xi \in \mathbb{R}^n : |\xi| \leq R \}$, so we can find another positive constant $\delta$ such that $$\frac{\operatorname{Re}\sigma + \kappa\,(\Lambda(\xi))^{2m-\rho}}{(\Lambda(\xi))^{2m}} \geq \delta, \hspace{0.3cm} |\xi| \leq R.$$
		So, the lemma is proved if we let $\eta = min (C, \delta)$.
	\end{pff}\\\\
	\textbf{Proof of Theorem \ref{Gardings theorem}} Let $T_{\tau} = J_{m}\,T_{\sigma}\,J_{m}$, where $J_m = T_{\sigma_{m}}$ and $\sigma_m(\xi) = (\Lambda(\xi))^{-m}$. Then, using the asymptotic expansion for the product of two pseudo-differential operators in Theorem 1.2 in \cite{MWMelliptic},
	$$T_{\sigma}\,J_m = T_{\tau_1},$$
	where
	\begin{eqnarray}\label{Equation with m-order}
		\tau_{1} - (\Lambda(.))^{-m} \sigma \in  M_{\rho,\Lambda}^{m-\rho}.
	\end{eqnarray}
	Similarly,
	$$T_{\tau}=J_{m}\,T_{\tau_{1}}$$
	and
	\begin{eqnarray}\label{Equation with 0-order}
		\tau - (\Lambda(.))^{-m} \tau_1 \in  M_{\rho,\Lambda}^{-\rho}.
	\end{eqnarray}
	Multiplying Equation \eqref{Equation with m-order} by $(\Lambda(.))^{-m}$ and adding the result to Equation \eqref{Equation with 0-order}, we get
	$$\tau - (\Lambda(.))^{-2m} \sigma \in  M_{\rho,\Lambda}^{-\rho}.$$
	Therefore
	$$
	\tau = (\Lambda(.))^{-2m} \sigma + r,
	$$
	where $r \in M_{\rho,\Lambda}^{-\rho}$. So, by Lemma \ref{strong elliptic inequality},
	$$
	\operatorname{Re} \tau = (\Lambda(.))^{-2 m}\,\operatorname{Re} \sigma + \operatorname{Re} r \geq \eta - \kappa(\Lambda(.))^{-\rho} + \operatorname{Re} r \geq \eta - \kappa^{\prime}\,(\Lambda(.))^{-\rho},$$
	where $\kappa^{\prime}$ is another positive constant. Therefore $\tau$ satisfies the conclusion of Lemma \ref{strong elliptic inequality} with $m=0$. Let us suppose for a moment that Gårding's inequality is valid for $m=0$. Then we can find a positive constant $C^{\prime}$ and a positive constant $C_{s}$ for every real number $s \geq \frac{\rho}{2}$ such that
	\begin{align*}
		\operatorname{Re}\left(T_{\sigma} \varphi, \varphi\right) &=\operatorname{Re}\left(J_{-m} T_{\tau} J_{-m} \varphi, \varphi\right) \\
		&=\operatorname{Re}\left(T_{\tau} J_{-m} \varphi, J_{-m} \varphi\right) \\
		& \geq C^{\prime}\,\|J_{-m} \varphi\|_{0,2,\Lambda}^{2} - C_{s}\,\|J_{-m} \varphi\|_{-s,2,\Lambda}^{2} \\
		&=C^{\prime}\,\|\varphi\|_{m, 2,\Lambda}^{2} - C_{s}\,\|\varphi\|_{m-s,2,\Lambda}^{2}
	\end{align*}
	for all $\varphi$ in $\mathcal{S}$. Now we need only to prove G\r{a}rding's inequality for $m=0$. Let $\sigma \in M_{\rho,\Lambda}^{0}$. Then, by Lemma \ref{strong elliptic inequality}, we can find positive constants $\eta$ and $\kappa$ such that
	$$
	\operatorname{Re} \sigma + \kappa\,(\Lambda(.))^{-\rho} \geq \eta.
	$$
	Let $F$ be a $C^{\infty}$ function on $\mathbb{C}$ such that
	$$
	F(z)=\sqrt{\frac{\eta}{2}+z},\hspace{0.3cm} z \in [0, \infty).
	$$
	Let $\tau$ be the function defined on $\mathbb{R}^{n} \times \mathbb{R}^{n}$ by
	$$
	\tau(x, \xi)=F\left(2\left(\operatorname{Re} \sigma(x, \xi)+\kappa\,(\Lambda(\xi))^{-\rho}-\eta\right)\right),\hspace{0.3cm} x, \xi \in \mathbb{R}^{n}.
	$$
	Then, by Lemma \ref{Compostion in M-class}, $\tau \in M_{\rho,\Lambda}^0$, and for all $x$ and $\xi$ in $\mathbb{R}^n$,
	\begin{align*}
		\tau(x, \xi) &= \sqrt{\frac{\eta}{2}+2 \operatorname{Re} \sigma(x, \xi)+2 \kappa\,(\Lambda(\xi))^{-\rho}-2 \eta} \\
		&=\sqrt{2 \operatorname{Re} \sigma(x, \xi)+2 \kappa(\Lambda(\xi))^{-\rho}-\frac{3}{2} \eta}.
	\end{align*}
	Using the asymptotic expansion for the formal adjoint of a pseudo- differential operator in Theorem 1.3 in \cite{MWMelliptic}, we get 
	$$
	T_{\tau}^{\ast}=T_{\tau^{\ast}},
	$$
	where $\tau^{\ast} \in M_{\rho,\Lambda}^0$ and $\tau-\tau^{\ast} \in M_{\rho,\Lambda}^{-\rho}$. Again, by using Theorem 1.2 in \cite{MWMelliptic}, we have
	$$
	T_{\tau}^{\ast} T_{\tau}=T_{\lambda},
	$$
	where
	$$
	\lambda-\tau^{\ast} \tau \in M_{\rho,\Lambda}^{-\rho}.
	$$
	If we let $r_{1}$ and $r_{1}^{\prime}$ in $ M_{\rho,\Lambda}^{-\rho}$ be such that
	$$
	\tau^{\ast}=\tau+r_{1}
	$$
	and
	$$
	\lambda=\tau^{\ast} \tau+r_{1}^{\prime}.
	$$
	Then
	$$
	\lambda=\left(\tau+r_{1}\right) \tau+r_{1}^{\prime}=2 \operatorname{Re} \sigma + 2 \kappa\,(\Lambda(.))^{-\rho}-\frac{3}{2} \eta + r_{2},
	$$
	where $$r_{2} = r_{1} \tau+r_{1}^{\prime} \in M_{\rho,\Lambda}^{-\rho}.$$
	So, if we let $r_{3}=2 \kappa\,(\Lambda(.))^{-\rho}+r_{2} \in M_{\rho,\Lambda}^{-\rho}$, then we get
	$$
	\lambda = 2 \operatorname{Re} \sigma-\frac{3}{2} \eta+r_{3}.
	$$
	But
	$$
	2 \operatorname{Re} \sigma=\sigma+\bar{\sigma}=\sigma+\sigma^{\ast}+r_{4}
	$$
	for some $r_{4}$ in $M_{\rho,\Lambda}^{-\rho}$. Therefore
	$$
	\lambda=\sigma+\sigma^{\ast}-\frac{3}{2} \eta+r_{5}
	$$
	for some $r_{5}$ in $M_{\rho,\Lambda}^{-\rho}$. Thus,
	$$
	\sigma+\sigma^{\ast}=\lambda+\frac{3}{2} \eta-r_{5}.
	$$
	Since
	$$
	\left(T_{\lambda} \varphi, \varphi\right)=\left(T_{\tau} \varphi, T_{\tau} \varphi\right) \geq 0, \hspace{0.3cm} \varphi \in \mathcal{S},
	$$
	it follows that
	\begin{align*}
		2 \operatorname{Re}\left(T_{\sigma}\, \varphi, \varphi\right) &=\left(T_{\sigma}\, \varphi, \varphi\right)+\left(T_{\sigma}^{\ast}\, \varphi, \varphi\right)=\left(T_{\sigma+\sigma^{\ast}} \,\varphi, \varphi\right) \\
		&=\left(T_{\lambda}\, \varphi, \varphi\right)+\frac{3}{2} \eta\,\|\varphi\|_{0,2,\Lambda}^{2}-\left(T_{r_{5}} \varphi, \varphi\right) \\
		& \geq \eta\,\|\varphi\|_{0,2,\Lambda}^{2}+\left\{\frac{\eta}{2}\,\|\varphi\|_{0,2,\Lambda}^{2}-\left\|T_{r_{5}}\, \varphi\right\|_{\frac{\rho}{2}, 2,\Lambda}\,\|\varphi\|_{-\frac{\rho}{2},2,\Lambda}\right\}, \hspace{0.3cm} \varphi \in \mathcal{S}.
	\end{align*}
	Then, by Theorem 1.6 in \cite{MWMelliptic} and using the fact that $r_5 \in M_{\rho,\Lambda}^{-\rho}$, we can find a positive constant $\nu$ such that
	$$2 \operatorname{Re}\left(T_{\sigma} \varphi, \varphi\right) \geq \eta\,\|\varphi\|_{0,2,\Lambda}^{2}+\left\{\frac{\eta}{2}\,\|\varphi\|_{0,2,\Lambda}^{2} - \nu\,\|\varphi\|_{-\frac{\rho}{2},2,\Lambda}^{2}\right\}, \hspace{0.3cm} \varphi \in \mathcal{S}.$$
	But
	$$
	\nu\|\varphi\|_{-\frac{\rho}{2}, 2,\Lambda}^{2}=\int_{\mathbb{R}^{n}} \nu\,(\Lambda(\xi))^{-\rho}\,|\hat{\varphi}(\xi)|^{2} d \xi=I+J,
	$$
	where
	$$
	I=\int_{\nu\,(\Lambda(\xi))^{-\rho} \leq \frac{\eta}{2}} \nu\,(\Lambda(\xi))^{-\rho}\,|\hat{\varphi}(\xi)|^{2}\, d \xi
	$$
	and
	$$
	J=\int_{\nu(\Lambda(\xi))^{-\rho} \geq \frac{\eta}{2}} \nu\,(\Lambda(\xi))^{-\rho}\,|\hat{\varphi}(\xi)|^{2}\, d \xi.
	$$
	Obviously,
	$$
	I \leq \frac{\eta}{2} \int_{\mathbb{R}^{n}}|\hat{\varphi}(\xi)|^{2} d \xi = \frac{\eta}{2}\,\|\varphi\|_{0,2,\Lambda}^{2}.
	$$
	To estimate $J$, we note that
	$$
	\nu\,(\Lambda(\xi))^{-\rho} \geq \frac{\eta}{2} \Rightarrow \Lambda(\xi) \leq \left(\frac{2 \nu}{\eta}\right)^{\frac{1}{\rho}}.
	$$
	So, for $\nu\,(\Lambda(\xi))^{-\rho} \geq \frac{\eta}{2}$, we get, for every real number $s \geq \frac{\rho}{2}$,
	\begin{align*}
		\nu\,(\Lambda(\xi))^{-\rho} & = \nu\,(\Lambda(\xi))^{2 s-\rho}\,(\Lambda(\xi))^{-2s} \\
		& \leq \nu\left(\frac{2 \nu}{\eta}\right)^{\frac{2 s-\rho}{\rho}}\,(\Lambda(\xi))^{-2s}.
	\end{align*}
	Thus for every real number $s \geq \frac{\rho}{2}$,
	$$
	J \leq \nu\left(\frac{2 \nu}{\eta}\right)^{\frac{2 s-\rho}{\rho}} \int_{\mathbb{R}^{n}}(\Lambda(\xi))^{-2 s}\,|\hat{\varphi}(\xi)|^{2}\, d \xi = C_{s}^{\prime}\,\|\varphi\|_{-s,2,\Lambda}^{2}.
	$$
	where $C_{s}^{\prime}=\nu\,\left(\frac{2 \nu}{\eta}\right)^{\frac{2 s-\rho}{\rho}}$. Therefore
	$$2\operatorname{Re}(T_{\sigma}\,\varphi,\varphi) \geq \eta\,\|\varphi\|_{0,2,\Lambda}^2 - C_{s}^{\prime}\,\|\varphi\|_{-s,2,\Lambda}^2,\hspace{0.3cm} \varphi \in \mathcal{S},$$
	and this completes the proof of the theorem.
	\begin{defi}
		Let $\sigma \in M_{\rho,\Lambda}^{m_1,m_2}$, where $m_1,m_2 \in \mathbb{R}$. Then $\sigma$ is said to be strongly SG $M$-elliptic if there exist positive constants $C$ and $R$ such that
		$$\operatorname{Re}(\sigma(x,\xi)) \geq C\Lambda(x)^{m_2}\Lambda(\xi)^{m_1},\hspace{.3cm} |x|^2 + |\xi|^2 \geq R^2.$$
	\end{defi}
	\begin{theorem}\textbf{(G\r{a}rding's Inequality For SG M-elliptic Operators)}\label{Gardings theorem for SG case}\\
		Let $\sigma \in M_{\rho,\Lambda}^{2m_1,2m_2}$, where $m_1,m_2 \in \mathbb{R}$, be strongly SG $M$-elliptic symbol. Then we can find  positive constants $C^{\prime}$ and $C_{s_1,s_2}$ for every real number $s_1  \leq \frac{\rho}{2}$ , $s_2 \geq \frac{\rho}{2}$ such that
		$$\operatorname{Re}(T_{\sigma}\,\varphi,\varphi) \geq C^{\prime}\,\|\varphi\|_{m_1,m_2,2,\Lambda}^2 - C_{s_1,s_2}\,\|\varphi\|_{m_1-s_1,m_2-s_2,2,\Lambda}^2,\hspace{0.3cm} \varphi \in \mathcal{S}.$$
	\end{theorem}
	To prove this theorem, we need the following two lemmas.
	\begin{lemma}\label{Compostion in SG M-class}
		Let $\sigma \in M_{\rho,\Lambda}^{0,0}$ and $F$ be a $C^{\infty}$- function on the complex plane $\mathbb{C}$. Then $F \circ \sigma \in M_{\rho,\Lambda}^{0,0}$.
	\end{lemma}
	Proof of the above lemma follows from the similar techniques as in Lemma \ref{Compostion in S-class} and Lemma \ref{Compostion in M-class}.
	\begin{lemma}\label{strong SG elliptic inequality}
		Let $\sigma \in M_{\rho,\Lambda}^{2m_1,2m_2}$, where $m_1,m_2 \in \mathbb{R}$, be strongly SG $M$-elliptic symbol. Then we can find two positive constants $\eta$ and $\kappa$ such that $$\operatorname{Re}(\sigma (x,\xi)) \geq \eta\,(\Lambda(x))^{2m_2}\,(\Lambda(\xi))^{2m_1} - \kappa\,(\Lambda(x))^{2m_2-\rho}\,(\Lambda(\xi))^{2m_1-\rho},\hspace{0.3cm} x,\xi \in \mathbb{R}^n.$$
	\end{lemma}
	Proof of the above lemma follows from the similar techniques as in Lemma \ref{strong elliptic inequality}.\\\\
	\textbf{Proof of Theorem \ref{Gardings theorem for SG case}} Let $T_{\tau} = J_{m_1,m_2}\,T_{\sigma}\,J_{m_1,m_2}$, where $J_{m_1,m_2} = T_{\sigma_{m_1,m_2}}$ and $\sigma_{m_1,m_2}(x,\xi) = (\Lambda(x))^{-m_2}\,(\Lambda(\xi))^{-m_1}$. Then, by Theorem \ref{product}, we have
	$$T_{\sigma}\,J_{m_1,m_2} = T_{\tau_1},$$
	where
	\begin{eqnarray}\label{Equation with m-order in SG}
		\tau_{1} - \sigma_{m_1,m_2} \sigma \in  M_{\rho,\Lambda}^{m_1-\rho,m_2-\rho}.
	\end{eqnarray}
	Similarly,
	$$T_{\tau}=J_{m_1,m_2}\,T_{\tau_{1}}$$
	and
	\begin{eqnarray}\label{Equation with 0-order in SG}
		\tau - \sigma_{m_1,m_2} \tau_1 \in  M_{\rho,\Lambda}^{-\rho,-\rho}.
	\end{eqnarray}
	Multiplying Equation \eqref{Equation with m-order in SG} by $\sigma_{m_1,m_2}$ and adding the result to Equation \eqref{Equation with 0-order in SG}, we get
	$$\tau - \sigma_{m_1,m_2}^2 \sigma \in  M_{\rho,\Lambda}^{-\rho,-\rho}.$$
	Therefore
	$$
	\tau =  \sigma_{m_1,m_2}^2 \sigma + r,
	$$
	where $r \in M_{\rho,\Lambda}^{-\rho,-\rho}$. So, by Lemma \ref{strong SG elliptic inequality},
	$$
	\operatorname{Re} \tau = \sigma_{m_1,m_2}^2\,\operatorname{Re} \sigma + \operatorname{Re} r \geq \eta - \kappa\,(\Lambda(x))^{-\rho}\,(\Lambda(\xi))^{-\rho}+ \operatorname{Re} r \geq \eta - \kappa^{\prime}\,(\Lambda(x))^{-\rho}\,(\Lambda(\xi))^{-\rho},
	$$
	where $\kappa^{\prime}$ is another positive constant. Therefore $\tau$ satisfies the conclusion of Lemma \ref{strong SG elliptic inequality} with $m_1=0$ and $m_2=0$. Let us suppose for a moment that Gårding's inequality is valid for $m_1=0$ and $m_2=0$. Then we can find a positive constant $C^{\prime}$ and a positive constant $C_{s_1,s_2}$ for every real number $s_1  \leq \frac{\rho}{2}$ , $s_2 \geq \frac{\rho}{2}$  such that
	\begin{align*}
		\operatorname{Re}\left(T_{\sigma} \varphi, \varphi\right) &=\operatorname{Re}\left(J_{-m_1,-m_2} T_{\tau}  J_{-m_1,-m_2}\,\varphi, \varphi\right) \\
		&=\operatorname{Re}\left(T_{\tau} J_{-m_1,-m_2}\,\varphi, J_{-m_1,-m_2}\,\varphi\right) \\
		& \geq C^{\prime}\,\|J_{-m_1,-m_2}\,\varphi\|_{0,0,2,\Lambda}^{2} - C_{s_1,s_2}\,\|J_{-m_1,-m_2}\,\varphi\|_{-s_1,-s_2,2,\Lambda}^{2} \\
		&=C^{\prime}\,\|\varphi\|_{m_1,m_2, 2,\Lambda}^{2} - C_{s_1,s_2}\,\|\varphi\|_{m_1-s_1,m_2-s_2,2,\Lambda}^{2}
	\end{align*}
	for all $\varphi$ in $\mathcal{S}$. Now we need only to prove G\r{a}rding's inequality for $m_1=0$ and $m_2=0$. Let $\sigma \in M_{\rho,\Lambda}^{0,0}$. Then, by Lemma \ref{strong SG elliptic inequality}, we can find positive constants $\eta$ and $\kappa$ such that
	$$
	\operatorname{Re} \sigma + \kappa\,(\Lambda(x))^{-\rho}\,(\Lambda(\xi))^{-\rho} \geq \eta.
	$$
	Let $F$ be a $C^{\infty}$ function on $\mathbb{C}$ such that
	$$
	F(z)=\sqrt{\frac{\eta}{2}+z},\hspace{0.3cm} z \in [0, \infty).
	$$
	Let $\tau$ be the function defined on $\mathbb{R}^{n} \times \mathbb{R}^{n}$ by
	$$
	\tau(x, \xi)=F\left(2\left(\operatorname{Re} \sigma(x, \xi)+\kappa\,(\Lambda(x))^{-\rho}\,(\Lambda(\xi))^{-\rho}-\eta\right)\right),\hspace{0.3cm} x, \xi \in \mathbb{R}^{n}.
	$$
	Then, by Lemma \ref{Compostion in SG M-class}, $\tau \in M_{\rho,\Lambda}^{0,0}$, and for all $x$ and $\xi$ in $\mathbb{R}^n$,
	\begin{align*}
		\tau(x, \xi) &= \sqrt{\frac{\eta}{2}+2 \operatorname{Re} \sigma(x, \xi)+2 \kappa\,(\Lambda(x))^{-\rho}\,(\Lambda(\xi))^{-\rho}-2 \eta} \\
		&=\sqrt{2 \operatorname{Re} \sigma(x, \xi)+2 \kappa\,(\Lambda(x))^{-\rho}\,(\Lambda(\xi))^{-\rho}-\frac{3}{2} \eta}.
	\end{align*}
	Hence, by Theorem \ref{adjoint}, we get 
	$$
	T_{\tau}^{\ast}=T_{\tau^{\ast}},
	$$
	where $\tau^{\ast} \in M_{\rho,\Lambda}^{0,0}$ and $\tau-\tau^{\ast} \in M_{\rho,\Lambda}^{-\rho,-\rho}$. Again, by using Theorem \ref{product}, we have
	$$
	T_{\tau}^{\ast} T_{\tau}=T_{\lambda},
	$$
	where
	$$
	\lambda-\tau^{\ast} \tau \in M_{\rho,\Lambda}^{-\rho,-\rho}.
	$$
	If we let $r_{1}$ and $r_{1}^{\prime}$ in $ M_{\rho,\Lambda}^{-\rho,-\rho}$ be such that
	$$
	\tau^{\ast}=\tau+r_{1}
	$$
	and
	$$
	\lambda=\tau^{\ast} \tau+r_{1}^{\prime}.
	$$
	Then
	$$
	\lambda=\left(\tau+r_{1}\right) \tau+r_{1}^{\prime}=2 \operatorname{Re} \sigma + 2 \kappa\,(\Lambda(x))^{-\rho}\,(\Lambda(\xi))^{-\rho}-\frac{3}{2} \eta + r_{2},
	$$
	where $$r_{2} = r_{1} \tau+r_{1}^{\prime} \in M_{\rho,\Lambda}^{-\rho,-\rho}.$$
	So, if we let $r_{3}=2 \kappa\,(\Lambda(x))^{-\rho}\,(\Lambda(\xi))^{-\rho}+r_{2} \in M_{\rho,\Lambda}^{-\rho,-\rho}$, then we get
	$$
	\lambda = 2 \operatorname{Re} \sigma-\frac{3}{2} \eta+r_{3}.
	$$
	But
	$$
	2 \operatorname{Re} \sigma=\sigma+\bar{\sigma}=\sigma+\sigma^{\ast}+r_{4}
	$$
	for some $r_{4}$ in $M_{\rho,\Lambda}^{-\rho,-\rho}$. Therefore
	$$
	\lambda=\sigma+\sigma^{\ast}-\frac{3}{2} \eta+r_{5}
	$$
	for some $r_{5}$ in $M_{\rho,\Lambda}^{-\rho,-\rho}$. Thus,
	$$
	\sigma+\sigma^{\ast}=\lambda+\frac{3}{2} \eta-r_{5}.
	$$
	Since
	$$
	\left(T_{\lambda} \varphi, \varphi\right)=\left(T_{\tau} \varphi, T_{\tau} \varphi\right) \geq 0, \hspace{0.3cm} \varphi \in \mathcal{S},
	$$
	it follows that
	\begin{align*}
		2 \operatorname{Re}\left(T_{\sigma}\, \varphi, \varphi\right) &=\left(T_{\sigma}\, \varphi, \varphi\right)+\left(T_{\sigma}^{\ast}\, \varphi, \varphi\right)=\left(T_{\sigma+\sigma^{\ast}} \,\varphi, \varphi\right) \\
		&=\left(T_{\lambda}\, \varphi, \varphi\right)+\frac{3}{2} \eta\,\|\varphi\|_{0,0,2,\Lambda}^{2}-\left(T_{r_{5}} \varphi, \varphi\right) \\
		& \geq \eta\,\|\varphi\|_{0,0,2,\Lambda}^{2}+\left\{\frac{\eta}{2}\,\|\varphi\|_{0,0,2,\Lambda}^{2}-\left\|T_{r_{5}}\, \varphi\right\|_{\frac{\rho}{2},\frac{\rho}{2}, 2,\Lambda}\,\|\varphi\|_{-\frac{\rho}{2},-\frac{\rho}{2},2,\Lambda}\right\}, \hspace{0.3cm} \varphi \in \mathcal{S}.
	\end{align*}
	Then, by Theorem \ref{bdd} and using the fact that $r_5 \in M_{\rho,\Lambda}^{-\rho,-\rho}$, we can find a positive constant $\nu$ such that
	$$2 \operatorname{Re}\left(T_{\sigma} \varphi, \varphi\right) \geq \eta\,\|\varphi\|_{0,0,2,\Lambda}^{2}+\left\{\frac{\eta}{2}\,\|\varphi\|_{0,0,2,\Lambda}^{2} - \nu\,\|\varphi\|_{-\frac{\rho}{2},-\frac{\rho}{2},2,\Lambda}^{2}\right\}, \hspace{0.3cm} \varphi \in \mathcal{S}.$$
	Since $\frac{\rho}{2} > 0$, so by Theorem \ref{contain-thm}, there exists a positive constant $C_0$ such that 
	$$|\varphi\|_{-\frac{\rho}{2},0,2,\Lambda}^2 \leq C_0\,\|\varphi\|_{0,0,2,\Lambda}^2, \hspace{0.3cm} \varphi \in \mathcal{S}.$$
	Now
	$$
	\nu\|\varphi\|_{-\frac{\rho}{2},-\frac{\rho}{2}, 2,\Lambda}^{2}=\int_{\mathbb{R}^{n}} \nu\,(\Lambda(x))^{-\rho}\,|(T_{\sigma_{\frac{\rho}{2}}}\varphi)(x)|^{2} dx=I+J,
	$$
	where
	$$
	I=\int_{\nu\,(\Lambda(x))^{-\rho} \leq \frac{\eta}{2C_0}} \nu\,(\Lambda(x))^{-\rho}\,|(T_{\sigma_{\frac{\rho}{2}}}\varphi)(x)|^{2} dx
	$$
	and
	$$
	J=\int_{\nu(\Lambda(x))^{-\rho} \geq \frac{\eta}{2C_0}} \nu\,(\Lambda(x))^{-\rho}\,|(T_{\sigma_{\frac{\rho}{2}}}\varphi)(x)|^{2} dx.
	$$
	Obviously,
	$$
	I \leq \frac{\eta}{2C_0} \int_{\mathbb{R}^{n}}|(T_{\sigma_{\frac{\rho}{2}}}\varphi)(x)|^{2}  dx = \frac{\eta}{2C_0}\,\|\varphi\|_{-\frac{\rho}{2},0,2,\Lambda}^2 \leq \frac{\eta}{2}\,\|\varphi\|_{0,0,2,\Lambda}^{2}.
	$$
	To estimate $J$, we note that
	$$
	\nu\,(\Lambda(x))^{-\rho} \geq \frac{\eta}{2C_0} \Rightarrow \Lambda(x) \leq \left(\frac{2 \nu C_0}{\eta}\right)^{\frac{1}{\rho}}.
	$$
	So, for $\nu\,(\Lambda(x)^{-\rho} \geq \frac{\eta}{2C_0}$, we get, for every real number $s_2 \geq \frac{\rho}{2}$,
	\begin{align*}
		\nu\,(\Lambda(x))^{-\rho} & = \nu\,(\Lambda(x))^{2 s_2-\rho}\,(\Lambda(x))^{-2s_2} \\
		& \leq \nu \left(\frac{2 \nu C_0}{\eta}\right)^{\frac{2s_2-\rho}{\rho}}\,(\Lambda(x))^{-2s_2}.
	\end{align*}
	Thus for every real number $s_2 \geq \frac{\rho}{2}$,
	$$
	J \leq \nu \left(\frac{2 \nu C_0}{\eta}\right)^{\frac{2s_2-\rho}{\rho}} \int_{\mathbb{R}^{n}}(\Lambda(x))^{-2 s_2}\,|(T_{\sigma_{\frac{\rho}{2}}}\varphi)(x)|^{2} dx = C_{s_2}^{\prime}\,\|\varphi\|_{-\frac{\rho}{2},-s_2,2,\Lambda}^{2},
	$$
	where $C_{s_2}^{\prime}=\nu \left(\frac{2 \nu C_0}{\eta}\right)^{\frac{2s_2-\rho}{\rho}}$. Therefore
	$$2\operatorname{Re}(T_{\sigma}\,\varphi,\varphi) \geq \eta\,\|\varphi\|_{0,0,2,\Lambda}^2 - C_{s_2}^{\prime}\,\|\varphi\|_{-\frac{\rho}{2},-s_2,2,\Lambda}^2,\hspace{0.3cm} \varphi \in \mathcal{S}.$$
	Since $\frac{\rho}{2} \geq s_1$, so by Theorem \ref{contain-thm}, there exists a positive constant $C_{s_1}^{\prime}$ such that 
	$$|\varphi\|_{-\frac{\rho}{2},-s_2,2,\Lambda}^2 \leq C_{s_1}^{\prime}\,\|\varphi\|_{-s_1,-s_2,2,\Lambda}^2, \hspace{0.3cm} \varphi \in \mathcal{S}.$$
	Thus for every real number $s_1  \leq \frac{\rho}{2}$ and $s_2 \geq \frac{\rho}{2}$, we have
	$$2\operatorname{Re}(T_{\sigma}\,\varphi,\varphi) \geq \eta\,\|\varphi\|_{0,0,2,\Lambda}^2 - C_{s_1,s_2}\,\|\varphi\|_{-\frac{\rho}{2},-s_2,2,\Lambda}^2,\hspace{0.3cm} \varphi \in \mathcal{S},$$
	where $C_{s_1,s_2} = C_{s_1}^{\prime}\,C_{s_2}^{\prime}$
	and this completes the proof of the theorem.
	

\begin{thebibliography}{99}
		\bibitem{Al2}M. Alimohammady and M. K. Kalleji. Spectral theory of a hybrid class of pseudo-differential operators,
		Complex Variables and Elliptic Equations, Volume 59, 12 2014. \vspace{10pt}
		\bibitem{AM} M. Alimohammady and M. K. Kalleji. Fredholmness property of M-elliptic pseudo-differential
		operator under change variable in its symbol. J. Pseudo-Differ. Oper. Appl., 4(3):371–392, 2013.\vspace{10pt}
		\bibitem{Bo} P. Boggiatto, E. Buzano and L. Rodino, Global Hypoellipticity and Spectral Theory,
		Akademie-Verlag, 1996.\vspace{10pt}
		\bibitem{AL} M. Cappiello and L. Rodino. SG-pseudodifferential operators and Gelfand-Shilov spaces. Rocky
		Mountain J. Math., 36(4):1117–1148, 2006.\vspace{10pt}
		\bibitem{AD} A. Dasgupta. Ellipticity of Fredholm pseudo-differential operators on  $L^p(\mathbb{R}^n)$. In New developments
		in pseudo-differential operators, volume 189 of Oper. Theory Adv. Appl., pages 107–116.
		Birkh\"{a}user, Basel, 2009.\vspace{10pt}
		\bibitem{AD&MW} A. Dasgupta and M. W. Wong. Spectral theory of SG pseudo-differential operators on  $L^p(\mathbb{R}^n)$.
		Studia Math., 187(2):185–197, 2008.\vspace{10pt}
		\bibitem{Yu} Yu. V. Egorov and B.-W. Schulze, Pseudo-Differential Operators, Singularities, Applications, Birkh\"auser, 1997.\vspace{10pt}
		\bibitem{GM} G. Garello and A. Morando. A class of $L^p$ bounded pseudodifferential operators. In Progress in
		analysis, Vol. I, II (Berlin, 2001), pages 689–696. World Sci. Publ., River Edge, NJ, 2003.\vspace{10pt}
		\bibitem{GM1} G. Garello and A. Morando, Lp-bounded pseudo-differential opearors and regularity
		for multi-quasi elliptic equations, Integr. Equ. Oper. Theory 51, 501-517 (2005).\vspace{10pt} 
		\bibitem{GA} G. Garello and A. Morando. m-microlocal elliptic pseudodifferential operators acting on $L_{loc}^p(\Omega)$.
		Math. Nachr., 289(14-15):1820–1837, 2016.\vspace{10pt}
		\bibitem{VVG} V. V. Gru\v{s}in. Pseudodifferential operators in $\mathbb{R}^n$ with bounded symbols. Funkcional. Anal. i
		Prilo\v{z}en, 4(3):37–50, 1970.\vspace{10pt}
		\bibitem{Hor}  H\"ormander, L., The analysis of linear partial differential operators, Vol. III. 
		Springer-Verlag, Berlin Heidelberg NewYork Tokyo, 1985. Pseudo-differential operators.\vspace{10pt}
		\bibitem{Ni} F. Nicola, K-theory of SG-pseudo-differential algebras, Proc. Amer. Math. Soc. 131
		(2003), 2841–2848.\vspace{10pt}
		\bibitem{NiRo} F. Nicola and L. Rodino, SG pseudo-differential operators and weak hyperbolicity,
		Pliska Stud. Math. Bulgar. 15 (2003), 5–20.\vspace{10pt}
		\bibitem{Ta} M. E. Taylor, Pseudo-differential operators, Princton university press, Princton
		(1981).\vspace{10pt}
		\bibitem{MWfredholm} M. W. Wong. Fredholm pseudodifferential operators on weighted Sobolev spaces. Ark. Mat.,
		21(2):271–282, 1983.\vspace{10pt}
		\bibitem{MWspectral} M. W. Wong. Spectral theory of pseudo-differential operators. Adv. in Appl. Math., 15(4):437–
		451, 1994.\vspace{10pt}
		\bibitem{MWMelliptic} M. W. Wong. M-elliptic pseudo-differential operators on $L^p(\mathbb{R}^n)$. Math. Nachr., 279(3):319–
		326, 2006.\vspace{10pt}
		\bibitem{MWbook} M. W. Wong. An introduction to pseudo-differential operators, volume 6 of Series on Analysis,
		Applications and Computation. World Scientific Publishing Co. Pte. Ltd., Hackensack, NJ,
		third edition, 2014.
	\end{thebibliography}
\end{document}